%%%%%%%%%%%%%%%%%%%%%%%%%%%%%%%%%%%%%%%%%%%%%%%%%%%%%%%%%%%%%%%%%%%%%%%%%%
%                                                                        % 
%        by J. Bouttier, P. Di Francesco and  E. Guitter                 %
%                TEX file, using lanlmac.tex macros                      %
%                                                                        %
%                                                                        %
%%%%%%%%%%%%%%%%%%%%%%%%%%%%%%%%%%%%%%%%%%%%%%%%%%%%%%%%%%%%%%%%%%%%%%%%%%
\input lanlmac
\def\href#1#2{{#2}}

\input epsf.tex

\overfullrule=0mm

\newcount\figno
\figno=0
\def\fig#1#2#3{
\par\begingroup\parindent=0pt\leftskip=1cm\rightskip=1cm\parindent=0pt
\baselineskip=11pt
\global\advance\figno by 1
\midinsert
\epsfxsize=#3
\centerline{\epsfbox{#2}}
\vskip 12pt
{\bf Fig.\ \the\figno:} #1\par
\endinsert\endgroup\par
}
\def\figlabel#1{\xdef#1{\the\figno}}
\def\encadremath#1{\vbox{\hrule\hbox{\vrule\kern8pt\vbox{\kern8pt
\hbox{$\displaystyle #1$}\kern8pt}
\kern8pt\vrule}\hrule}}

%Macros 
%%%%%%%%%%%%%%%%%%%%%%%%%%%%%%%%%%%%%%%%%%%%%%%%%%%%%%%%%%%%%%%%%

\def\IR{\relax{\rm I\kern-.18em R}}
\font\cmss=cmss10 \font\cmsss=cmss10 at 7pt

\font\cmss=cmss10 \font\cmsss=cmss10 at 7pt
\def\IZ{\relax\ifmmode\mathchoice
{\hbox{\cmss Z\kern-.4em Z}}{\hbox{\cmss Z\kern-.4em Z}}
{\lower.9pt\hbox{\cmsss Z\kern-.4em Z}}
{\lower1.2pt\hbox{\cmsss Z\kern-.4em Z}}\else{\cmss Z\kern-.4em Z}\fi}
\def\IN{\relax{\rm I\kern-.18em N}}
\def\circbullet{{\bigcirc \kern-.75em \bullet \kern .3em}}
\def\circcirc{{\bigcirc \kern-.75em \circ \kern.3em}}
\def\smallcircbullet{{{\scriptscriptstyle{\bigcirc \kern-.5em \bullet}} \kern .2em}}
\def\smallcirccirc{\relax{{\scriptscriptstyle{\bigcirc \kern-.5em \circ}} \kern .2em}}
\def\b{\circ}
\def\n{\bullet}

\def\gbbbb{\Gamma_4^{\hbox{$\scriptstyle \b \b$}\kern -8.2pt
\raise -4pt \hbox{$\scriptstyle \b \b$}}}
\def\gnnnn{\Gamma_4^{\hbox{$\scriptstyle \n \n$}\kern -8.2pt  
\raise -4pt \hbox{$\scriptstyle \n \n$}}}
\def\gnnnnnn{\Gamma_6^{\hbox{$\scriptstyle \n \n \n$}\kern -12.3pt
\raise -4pt \hbox{$\scriptstyle \n \n \n$}}}
\def\gbbbbbb{\Gamma_6^{\hbox{$\scriptstyle \b \b \b$}\kern -12.3pt
\raise -4pt \hbox{$\scriptstyle \b \b \b$}}}
\def\gbbbbc{\Gamma_{4\, c}^{\hbox{$\scriptstyle \b \b$}\kern -8.2pt
\raise -4pt \hbox{$\scriptstyle \b \b$}}}
\def\gnnnnc{\Gamma_{4\, c}^{\hbox{$\scriptstyle \n \n$}\kern -8.2pt
\raise -4pt \hbox{$\scriptstyle \n \n$}}}
\def\Rbud#1{{\cal R}_{#1}^{-\kern-1.5pt\blacktriangleright}}
\def\Rleaf#1{{\cal R}_{#1}^{-\kern-1.5pt\vartriangleright}}
\def\Rbudb#1{{\cal R}_{#1}^{\circ\kern-1.5pt-\kern-1.5pt\blacktriangleright}}
\def\Rleafb#1{{\cal R}_{#1}^{\circ\kern-1.5pt-\kern-1.5pt\vartriangleright}}
\def\Rbudn#1{{\cal R}_{#1}^{\bullet\kern-1.5pt-\kern-1.5pt\blacktriangleright}}
\def\Rleafn#1{{\cal R}_{#1}^{\bullet\kern-1.5pt-\kern-1.5pt\vartriangleright}}
\def\Wleaf#1{{\cal W}_{#1}^{-\kern-1.5pt\vartriangleright}}
\def\Cleaf{{\cal C}^{-\kern-1.5pt\vartriangleright}}
\def\Cbud{{\cal C}^{-\kern-1.5pt\blacktriangleright}}
\def\Crleaf{{\cal C}^{-\kern-1.5pt\circledR}}

%%%%%%%%%%%%%%%%%%%%%%%%%%%%%%%%%%%%%%%%%%%%%%%%%%%%%%%%%%%%%%%%%
%%%%%%%%%%%%%%%%%%%%%%%%%%%%%%%%%%%%%%%%%%%%%%%%%%%%%%%%%%%%%%%%%%%%%

\magnification=\magstep1
\baselineskip=12pt
\hsize=6.3truein
\vsize=8.7truein
 at 8truept
 at 8truept
 at 10truept

%%%%%%%%%%%%%%%%%%%%%%%%%%%%%%%%%%%%%%%%%%%%%%%%%%%%%%%%%%%%%%%%%%%%%%%%
\font\bigrm=cmr12 at 14pt
\centerline{\bigrm Combinatorics of bicubic maps with hard particles}

\bigskip\bigskip

\centerline{J. Bouttier${}^{1,2}$, P. Di Francesco${}^1$ and E. Guitter${}^1$}
  \medskip
  \centerline{\it ${}^1$Service de Physique Th\'eorique, CEA/DSM/SPhT}
  \centerline{\it Unit\'e de recherche associ\'ee au CNRS}
  \centerline{\it CEA/Saclay}
  \centerline{\it 91191 Gif sur Yvette Cedex, France}
  \smallskip
  \centerline{\it ${}^2$Instituut voor Theoretische Fysica}
  \centerline{\it Valckeniertstraat 65}
  \centerline{\it 1018 XE Amsterdam}
  \centerline{\it The Netherlands}
  \medskip
\centerline{\tt bouttier@spht.saclay.cea.fr}
\centerline{\tt philippe@spht.saclay.cea.fr}
\centerline{\tt guitter@spht.saclay.cea.fr}

  \bigskip
  \bigskip
  \bigskip

%  \centerline{\footrm Submitted: Jan 1, 2004;  Accepted: Jan 2, 2004;
%     Published: Jan 3, 2004}
%     \centerline{\footrm Mathematics Subject Classifications: Primary 05C30; Secondary 05A15,
%     05C05, 05C12, 68R05}

     \bigskip\bigskip

     \centerline{\bf Abstract}
     \smallskip
     {\narrower\noindent
We present a purely combinatorial solution of the problem of enumerating planar bicubic
maps with hard particles. This is done by use of a bijection with a particular class of blossom 
trees with particles, obtained by an appropriate cutting of the maps. Although these trees have 
no simple local characterization, we prove that their enumeration may be performed upon introducing 
a larger class of ``admissible" trees with possibly doubly-occupied edges and summing them with 
appropriate signed weights. The proof relies on an extension of the cutting procedure 
allowing for the presence on the maps of special non-sectile edges. The admissible trees are 
characterized by simple local rules, allowing eventually for an exact enumeration of planar bicubic maps 
with hard particles. We also discuss generalizations for maps with particles subject to
more general exclusion rules and show how to re-derive the enumeration of quartic maps with Ising spins
in the present framework of admissible trees. We finally comment on a possible interpretation
in terms of branching  processes.
\par}

   \bigskip
   \bigskip
   \bigskip

%AMS Subject Classification (2000): Primary 05C30; Secondary 05A15, 
%05C05, 05C12, 68R05

%references
\nref\TUT{W. Tutte, {\it A Census of Planar Maps}, Canad. J. of Math. {\bf 15} (1963) 249-271;
{\it A Census of planar triangulations}, Canad. J. of Math. {\bf 14} (1962) 21-38;
{\it A Census of Hamiltonian polygons}, Canad. J. of Math. {\bf 14} (1962) 402-417;
{\it A Census of slicings}, Canad. J. of Math. {\bf 14} (1962) 708-722.}
\nref\BIPZ{E. Br\'ezin, C. Itzykson, G. Parisi and J.-B. Zuber, {\it Planar
Diagrams}, Comm. Math. Phys. {\bf 59} (1978) 35-51.}
\nref\DGZ{P. Di Francesco, P. Ginsparg and J. Zinn--Justin, {\it 2D Gravity and Random Matrices},
Physics Reports {\bf 254} (1995) 1-131.}
\nref\Bessis{D. Bessis, {\it A new method in the combinatorics of the topological expansion},
Comm. Math. Phys. {\bf 69} (1979) 147-163; D. Bessis, C. Itzykson and J.-B. Zuber, {\it Quantum 
field theory techniques in graphical enumeration}, Adv. Appl. Math {\bf 1} (1980) 109-157.}
\nref\Amb{J. Ambj\o rn, B. Durhuus and T. Jonsson, {\it Quantum geometry. A statistical field
theory approach}, Cambridge Monogr. Math. Phys. {\bf 1} (1997).}
\nref\ACKM{J. Ambj\o rn, L. Chekhov, C. Kristjansen and Y. Makeenko, {\it Matrix model
calculations beyond the spherical limit}, Nucl. Phys. B {\bf 404} (1993) 127-172. [Erratum-ibid
{\bf B 449} (1995) 681].}
\nref\CORV{R. Cori and B. Vauquelin, {\it Planar maps are well labeled trees},
Canad. J. Math. {\bf 33 (5)} (1981) 1023-1042.}
\nref\ARQUES{D. Arqu\`es, {\it Les hypercartes planaires sont des arbres 
tr\`es bien \'etiquet\'es}, Discr. Math. {\bf 58(1)} (1986) 11-24.}  
\nref\SCHth{G. Schaeffer, {\it Conjugaison d'arbres
et cartes combinatoires al\'eatoires}, PhD Thesis, Universit\'e 
Bordeaux I (1998).}
\nref\SCH{G. Schaeffer, {\it Bijective census and random
generation of Eulerian planar maps}, Elec.
Jour. of Combinatorics Vol. {\bf 4} (1997) R20.}
\nref\CONST{M. Bousquet-M\'elou and G. Schaeffer,
{\it Enumeration of planar constellations}, Adv. in Applied Math.,
{\bf 24} (2000) 337-368.}
\nref\COL{J. Bouttier, P. Di Francesco and E. Guitter, {\it Counting colored random
triangulations}, Nucl. Phys. {\bf B641[FS]} (2002) 519-532, arXiv:cond-mat/0206452.}
\nref\CENSUS{J. Bouttier, P. Di Francesco and E. Guitter, {\it Census of planar
maps: from the one-matrix model solution to a combinatorial proof},
Nucl. Phys. {\bf B645}[PM] (2002) 477-499, arXiv:cond-mat/0207682.}
\nref\BMS{M. Bousquet-M\'elou and G. Schaeffer,{\it The degree distribution
in bipartite planar maps: application to the Ising model},
arXiv:math.CO/0211070.}
\nref\CHP{J. Bouttier, P. Di Francesco and E. Guitter, {\it Combinatorics of 
hard particles on planar graphs}, Nucl. Phys. {\bf B655}[FS] (2003) 313-341, 
arXiv:cond-mat/0211168.}
\nref\KKMW{H. Kawai, N. Kawamoto, T. Mogami and Y. Watabiki, {\it Transfer matrix 
formalism for two-dimensional quantum gravity and fractal structures of space-time},
Phys. Lett. B {\bf 306} (1993) 19-26.}
\nref\Wata{Y. Watabiki, {\it Construction of noncritical string field theory by 
transfer matrix formalism in dynamical triangulation}, Nucl. Phys. B {\bf 441} 
(1995) 119-166.}
\nref\AW{J. Ambj\o rn and Y. Watabiki, {\it Scaling in quantum gravity}, 
Nucl. Phys. B {\bf 445} (1995) 129-144.}
\nref\AKW{J. Ambj\o rn, C. Kristjansen and Y. Watabiki, {\it The two-point 
function of $c = -2$ matter coupled to 2D quantum gravity}, 
Nucl.\ Phys.\ B {\bf 504} (1997) 555-578.}
\nref\CS{P. Chassaing and G. Schaeffer, {\it Random Planar Lattices and 
Integrated SuperBrownian Excursion}, 
Probability Theory and Related Fields {\bf 128(2)} (2004) 161-212, 
arXiv:math.CO/0205226.}
\nref\GEOD{J. Bouttier, P. Di Francesco and E. Guitter, {\it Geodesic
distance in planar graphs}, Nucl. Phys. {\bf B663}[FS] (2003) 535-567, 
arXiv:cond-mat/0303272.}
\nref\LALLER{J. Bouttier, P. Di Francesco and E. Guitter, {\it Random
trees between two walls: Exact partition function}, J. Phys. A: Math. Gen.
{\bf 36} (2003) 12349-12366, arXiv:cond-mat/0306602.} 
\nref\ONEWALL{J. Bouttier, P. Di Francesco and E. Guitter, {\it Statistics
of planar maps viewed from a vertex: a study via labeled trees},
Nucl. Phys. {\bf B675}[FS] (2003) 631-660, arXiv:cond-mat/0307606.}
\nref\CHASDUR{P. Chassaing and B. Durhuus, {\it Statistical Hausdorff 
dimension of labelled trees and quadrangulations}, 
arXiv:math.PR/0311532.}
\nref\DOU{see for instance M. Douglas, {\it The two-matrix model}, in 
{\it Random Surfaces and Quantum Gravity}, O. Alvarez,
E. Marinari and P. Windey eds., NATO ASI Series {B:} Physics Vol. {\bf 262} (1990).} 
\nref\GFRC{D. Gaunt and M. Fisher,
{\it Hard-Sphere Lattice Gases.I.Plane-Square Lattice}, J. Chem. Phys.
{\bf 43} (1965) 2840-2863; L. Runnels, L. Combs and J. Salvant, {\it Exact Finite
Methods of Lattice Statistics. II. Honeycomb-Lattice Gas of Hard Molecules},
J. Chem. Phys. {\bf 47} (1967) 4015-4020.}
\nref\BaxHH{R. J. Baxter, {\it Hard Hexagons: Exact Solution}, J. Phys. {\bf A 13}
(1980) L61-L70; R. J. Baxter and S.K. Tsang, {\it Entropy of Hard Hexagons},
J. Phys. {\bf A 13} (1980) 1023-1030; see also
R. J. Baxter, {\it Exactly Solved Models in Statistical Mechanics},
Academic Press, London (1984).}
\nref\BaxHS{R. J. Baxter, I. G. Enting and S.K. Tsang, {\it Hard Square Lattice Gas},
J. Stat. Phys. {\bf 22} (1980) 465-489; G. Kamieniarz and H. Bl\"ote, {\it The Non-interacting
Hard-square Lattice Gas: Ising Universality}, J. Phys. A Math. Gen.
{\bf 26} (1993) 6679-6689.}
\nref\CRI{J. Bouttier, P. Di Francesco and E. Guitter, {\it Critical and tricritical
hard objects on bicolourable random lattices: exact solutions}, J. Phys. A: Math. Gen.
{\bf 35} (2002) 3821-3854, arXiv:cond-mat/0201213.}
\nref\HARM{P. Di Francesco, {\it Geometrically Constrained Statistical Models on
Fixed and Random Lattices: From Hard Squares to Meanders}, arXiv:cond-mat/0211591.}
\nref\ISING{D. Boulatov and V. Kazakov, {\it The Ising model
on a random planar lattice: the structure of the phase
transition and the exact critical exponents}, Phys. Lett. {\bf B186} (1987) 379-384.}
\nref\GALWA{S. Karlin and H. Taylor, {\it A first course in
stochastic processes}, Academic Press, New-York (1975).}
\nref\MOB{J. Bouttier, P. Di Francesco and E. Guitter, {\it Planar
maps as labeled mobiles}, Elec. Jour. of Combinatorics Vol. {\bf 11} (2004) R69,
arXiv:math.CO/0405099.}
\nref\Aldous{D. Aldous, {\it Tree-Based Models for Random Distribution of Mass}, 
J. Stat. Phys. {\bf 73} (1993) 625-641.}

%text
\newsec{Introduction}

Enumeration of planar maps has been a constant subject of interest for both mathematicians 
and physicists. Many exact enumeration formulae have been obtained over the years for various 
types of maps, say maps with fixed vertex and/or face valences, maps with particular color 
prescriptions, or even maps carrying additional degrees of freedom such as spins or particles. Since 
the original combinatorial approach of Tutte \TUT, two general methods of
enumeration have been developed. The first one uses {\it matrix integrals} [\xref\BIPZ,\xref\DGZ] 
of the form $\int {\rm Exp}[-N{\rm Tr} U]$ over possibly several $N\times N$ Hermitian matrices and
with an action $U$ designed so that the diagrammatic expansion precisely reproduces the maps 
to be enumerated. 
For a large class of maps, a direct computation of the associated matrix integral 
may be carried out exactly, thus providing explicit enumeration formulae. 
These calculations rely on various techniques: saddle point methods \BIPZ, orthogonal
polynomial methods \Bessis, as well as the so-called loop equation method [\xref\Amb,\xref\ACKM], 
very similar in spirit to Tutte's original work. Beyond matrix integrals, another,
more recently developed technique relies on the existence of {\it bijections} between planar maps and 
decorated trees, reducing the problem to that, much simpler, of counting trees [\xref\CORV-\xref\CHP]. 
In most cases,
the trees at hand are easily constructed recursively and recursion formulae may therefore be
derived for the generating functions of the maps. The main interest of this technique is clearly
its conceptual simplicity, as it uses only elementary combinatorics. Another advantage is that it 
also gives access to refined properties of the maps. For instance, the statistics of distances 
between points on the maps, already addressed in Refs.[\xref\KKMW-\xref\AKW] by a specific
combinatorial treatment, emerges naturally from the bijective formalism [\xref\CS-\xref\CHASDUR].

The matrix integral and bijective techniques above are {\it a priori} very different in nature but, 
remarkably enough, 
in all cases which have been solved by both methods, the precise form of the solution
involves very similar systems of algebraic equations. That algebraic equations appear
as recursive equations for generating functions of trees should not come as a surprise. What
is more surprising is that the matrix integral solution comes precisely in the {\it same} form, 
allowing to interpret the somewhat abstract auxiliary functions of this matrix integral 
solution as generating functions for trees. Even if we have no direct understanding of this 
phenomenon, this simple observation allows us to use the matrix integral solution for a particular 
map enumeration problem as a tool to {\it guess} how to characterize the ensemble of trees in 
bijection with these maps. This already proved very useful for various classes of maps and 
the present paper will show another example where applying such a principle leads to the combinatorial, 
bijective solution of a quite involved problem, namely the combinatorics of bicubic maps
with hard particles.

Before we proceed to our study, let us recall more precisely to which families of maps
the methods above have been applied so far. We may classify the problems in the matrix integral
language according to the number of (Hermitian) matrices needed to reproduce the desired maps.
The simplest case of the so-called one-matrix integral over a matrix $M$ with action
\eqn\onemat{U(M)={ M^2\over 2}-\sum_{k\geq 1} g_k {M^k\over k}}
corresponds to maps with prescribed weights $g_k$ per $k$-valent vertex and was solved extensively 
in \BIPZ. Its tree formulation involves
so-called {\it blossom trees} \SCHth\ with ``charged" decorations, and with
simple charge characterizations \CENSUS.
The next case is that of the so-called two-matrix model involving an integral over
two matrices $B$ and $W$ interacting only via a term $BW$, namely with action
\eqn\twomat{U(B,W)=BW -\sum_{k\geq 1} g_k {B^k\over k}-\sum_{k\geq 1} {\tilde g}_k {W^k\over k}}
This model corresponds to
having {\it bipartite} maps, i.e.\ maps with two types of vertices, say black and white,
such that all edges connect vertices of opposite colors, and with prescribed weights $g_k$,
${\tilde g}_k$ according to the valence {\it and} the color of the vertices.
This matrix model was also solved extensively \DOU. Its
tree version involves again blossom trees with now slightly more involved 
charge characterizations [\xref\BMS,\xref\CHP]. Note that the one-matrix model may be viewed as a particular
case of the two-matrix model with $g_k=\delta_{k,2}$, in which case the matrix $B$ may be 
formally integrated out. In the map language, this simply states that arbitrary
maps are in bijection with bipartite maps whose black vertices all have valence $2$, as is
clear by inserting a bivalent black vertex in the middle of each edge. More generally,
we may always decide to suppress all (black or white) bivalent vertices of a bipartite map,
thus creating maps whose black and white vertices may be connected regardless of their color,
and which now carry edge weights depending on the relative colors of their adjacent vertices. 
The two-matrix model or its equivalent tree
formulation give for instance access to maps carrying configurations of the Ising model \BMS\ 
or maps with hard particles \CHP.

This exhausts the classes of maps which could be enumerated both via matrix integrals and
bijective methods. 
On the other hand, it is well known that multi-matrix integrals over, say $m$ matrices
$M_1\ldots M_m$ and with {\it chain interaction}, namely 
\eqn\multimat{U(M_1,\ldots,M_m)=\sum_{i=1}^{m-1} c_i M_i M_{i+1} -\sum_{i=1}^m{\sum_{k\geq 1} g^{(i)}_k 
{M_i^k\over k}}}
are completely calculable by means of bi-orthogonal polynomial techniques and give rise
to solutions expressible in terms of abstract auxiliary functions solutions of algebraic
equations \DGZ. These turn out to describe quite generally maps with particles subject to 
some {\it exclusion rules}. More precisely, the above $m$-matrix model describes maps with
particles occupying the vertices and with the edge exclusion constraint that  
the {\it total} number of particles at the ends of any edge does not exceed $m-1$.
In the case of even $m=2p$, the model alternatively describes bipartite maps
with particles such that the total number of particles at the ends of any edge now does not 
exceed \foot{The equivalence between the two formulations is obtained as follows: starting
from the bipartite formulation with the edge rule of at most $p-1$ particle occupancy, 
we single out bivalent black empty vertices and erase them, resulting in an effective modified edge rule
that amounts to having at most $2p-1$ particles per edge upon reinterpreting black $\ell$-occupied
vertices as white $p+\ell$-occupied ones.} $p-1$. For all these models, the algebraic nature
of the solution suggests that again a bijection with trees should exist.

The aim of this paper is precisely to extend the bijective techniques to maps with 
particles subject to edge exclusion rules. Very generally, the physical motivation
for considering such models is to construct simple solvable models describing gases 
of interacting particles. This was first achieved on regular two-dimensional lattices 
with the so-called {\it hard-core lattice gas}, providing the simplest possible 
model of a dense fluid which incorporates the excluded volume effect \GFRC. An exact 
solution for this model was obtained by Baxter in \BaxHH, in the case of hard particles 
living on the vertices of the triangular lattice (the so-called hard hexagon model). 
It displays a crystallization transition from a low density disordered phase to a high density
ordered (crystalline) phase in which particles tend to occupy preferentially
one of the three canonical sublattices of the (tripartite) triangular lattice. 
The transition was identified to be in the universality class of the critical 3-state Potts model, 
in agreement with the 3-fold symmetry of the crystalline groundstates. The square and
hexagonal lattice versions of the hard-core lattice gas, corresponding
respectively to hard squares and hard triangles, could not be solved exactly.
Still, numerical evidence points to a crystallization transition now in the class of
the critical Ising model \BaxHS, consistent with the twofold symmetry of the groundstates,
in which particles occupy one of the two canonical sublattices of these
bipartite lattices. This hunt for exact solutions led to consider the hard-core gas
model on random lattices, in the form of planar maps with hard particles.
It was recognized in \CRI\ that the crystallization transition is wiped out for the statistical
ensemble of arbitrary maps but that it is restored when restricting the class of maps to bipartite ones.
Simple examples are hard particles on planar bicubic or biquartic maps, which may be seen as the 
randomized, bipartite versions of the hard triangle and hard square models. These were solved exactly 
in \CRI\ and shown to display a crystallization transition in the class of the critical Ising model 
on random graphs. Indeed, the bipartite nature of the maps grants the existence of 
two crystalline maximally occupied groundstates, with particles lying on one of the two canonical
subgraphs of the maps. More generally, the above models
of particles on bipartite maps with the exclusion rule of at most $p-1$ particles per edge
allow to attain multi-critical transition points. When $p=3$ for instance, a tricritical
point in the class of the tricritical Ising model may be reached \CRI, separating a line
of Ising like transition points from a line of first order transition points between groundstates
of different symmetry.

In this paper, we concentrate mainly
on {\it planar bicubic maps with hard particles}, i.e.\ bipartite maps with black and white 
trivalent vertices and the exclusion rule that each edge has at most one adjacent particle.
In the graph theory language, planar bicubic maps correspond to planar bipartite 3-regular
(multi-)graphs and a configuration of hard particles in such a map corresponds to
an independent or stable set in the graph.
In the matrix language, this model is described by a four-matrix integral with action
\eqn\fourmat{U(M_1,\ldots,M_4)=M_1 M_2-M_2 M_3+M_3 M_4 -g \left({M_1^3\over 3}+{M_4^3\over 3}\right)\
-g\,z \left( {M_2^3\over 3}+{M_3^3\over 3}\right)}
where $g$ is a weight per vertex and $z$ a weight per particle. The explicit algebraic
solution of this matrix integral
was given in Ref.\CRI. In this paper, we show how to re-derive this solution in a purely
combinatorial way using bijections with suitable classes of blossom trees with particles.
As opposed to bipartite maps without particles (two-matrix model), a new feature of 
this construction is that we have to resort to a kind of inclusion-exclusion principle 
in which doubly-occupied edges are {\it a priori} allowed but are effectively subtracted 
by appropriate sign factors. 

The paper is organized as follows. In Sect.2, we recall the definition of
bicubic maps with hard particles and present the matrix model results of Ref.\CRI,
as well as the underlying physics. Sect.3 presents our main theorem, namely
the existence of a blossom tree formulation of the model in terms of so-called 
``admissible trees" whose enumeration with appropriate signed weights reproduce 
that of planar bicubic maps with hard particles. The definition of
blossom trees without particles is recalled in Sect.3.1. We then define in Sect.3.2
admissible trees as blossom trees with particles and with specific charge characterizations
and state our main theorem relating the generating function of bicubic maps
with particles to that of admissible trees with signed weights.
The characterization of admissible trees is used to list all their possible local environments
(Sect.3.3) and to derive a system of algebraic equations which
determines their desired signed generating function (Sect.3.4), 
reproducing the results of Ref.\CRI. The proof of the map/tree equivalence is given 
in Sects.4 and 5. It is based on
a general cutting procedure of bicubic maps into trees (Sect.4.1) adapted so
as to allow for ``special" {\it non-sectile} edges on the map (Sect.4.2).
This generalized cutting procedure is used in Sect.4.3 to construct a bijection
between trees and maps, both with special edges. It is used in Sect.5.1 in
the case of bicubic maps with hard particles to organize admissible trees into
equivalence classes represented by ``admissible maps", of which bicubic maps
with hard particles are only a subset. The characterization of these equivalence
classes (Sect.5.2) allows to show that, in the enumeration of admissible trees 
with appropriate signed weights, the contribution of an equivalence class is 
zero unless the admissible map satisfies the hard particle constraint (Sect.5.3).
This proves that the generating function for signed admissible trees matches
that of bicubic maps with hard particles.
Finally, we discuss in Sect.6 various generalizations of this map/tree equivalence
for the Ising model on quartic maps (Sect.6.1) or maps with generalized edge
exclusion rules (Sect.6.2). We also discuss a corollary of our construction, namely
the existence of an Ising-like crystallization transition for a statistical ensemble 
of trees with particles and comment on a possible branching process interpretation.
We gather a few concluding remarks in Sect.7.

\fig{An example of bicubic map with hard particles. The vertices are of four types: black-empty, 
white-empty, black-occupied and white-occupied, represented respectively by $\bullet$, $\circ$, 
$\circbullet$ and $\circcirc$. All vertices have valence three and adjacent vertices have different
colors. The hard particle constraint states that no two adjacent vertices are simultaneously
occupied. The map is rooted at the edge indicated by the arrow, originating from a
black-empty vertex and with the external face on its left.}{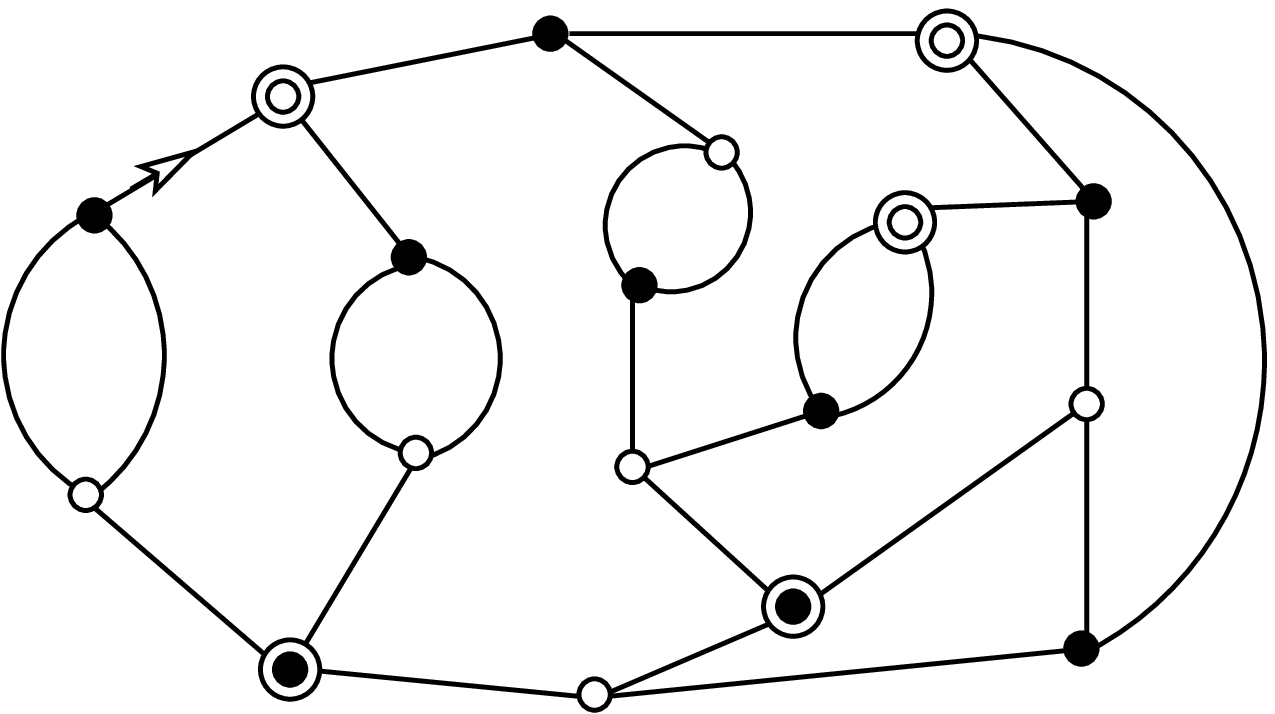}{8.cm}
\figlabel\bicub
\newsec{The problem of bicubic maps with hard particles and its matrix integral solution}
\subsec{Hard particles on bicubic maps}
As mentioned in the introduction, we will mainly study in this paper 
{\it bicubic maps with hard particles}, namely planar
maps with black and white trivalent vertices and edges connecting only vertices of
opposite color. The vertices may be occupied or not by particles, with the {\it hard
particle constraint} that no edge has two adjacent occupied vertices. The map is
assumed to contain at least one black-empty vertex and it is furthermore rooted,
i.e.\ carries a marked oriented edge pointing from one of its black-empty vertices (the root vertex).
An example of such a map is depicted in Fig.\bicub. In the planar representation, we always
chose for external face that on the left of the root edge.
We are interested in enumerating these maps with, say $n_\circ$ white-empty, $n_\bullet +1$
black-empty, $n_\smallcirccirc$ white-occupied and $n_\smallcircbullet$ black-occupied vertices.
For simplicity, we will concentrate here on the slightly easier task of computing the 
{\it generating function} $G_{\rm BMHP}(g,z)$ for these maps with, say a weight $g$ per empty 
or occupied, black or white vertex, and a weight $z$ per particle, namely
\eqn\genexpl{G_{\rm BMHP}(g,z)=\sum_{m,n} g^n z^m {\cal B}_{n,m}}
where ${\cal B}_{n,m}$ denotes the total number of planar bicubic maps rooted at a black-empty
vertex with a total number $n=n_\circ+n_\bullet+n_\smallcirccirc+n_\smallcircbullet +1$ of vertices 
and a total number $m=n_\smallcirccirc+n_\smallcircbullet$ of particles.

\subsec{Matrix integral solution}

The generating function $G_{\rm BMHP}(g,z)$ was obtained in Ref.\CRI\ by use of
a four-matrix integral with action \fourmat. 
The result involves four auxiliary functions $R$, $S$, $T$ and $V$
depending on $g$ and $z$ and solution of a closed algebraic system: 
\eqn\matrecur{\eqalign{
S&= 1 + 2 g\, R \cr T &= -g^3z\, V^4 \cr R &= g\, V^2+gz\, S^2 +2 g^2z\, T\cr
V&= S+2g^2z\, VS\cr}}
which determines them uniquely order by order in $g$ provided we take $V=1+{\cal O}(g)$.
The planar free energy $F_{\rm BMHP}(g,z)$ of the model, counting {\it unrooted} bicubic maps
with hard particles, and without the constraint of having at least one black-empty vertex,
is related to these functions via
\eqn\freeenerg{ \left({1\over g} {d \over dg} \right)^2 g^4 F_{\rm BMHP}(g,z)=4  {\rm Log}
V(g,z)}
The generating function $G_{\rm BMHP}(g,z)$ is related to the free energy via
\eqn\FtoG{G_{\rm BMHP}(g,z)={3\over 2}\left(g {d \over dg}-z {d \over dz} \right) F_{\rm BMHP}(g,z)}
where the derivatives amount to marking a vertex ($g\, d/dg$) which is not occupied
($-z\, d/dz$). The prefactor $1/2$ amounts to assuming that this vertex is black
and the factor $3$ stands for the $3$ choices of an adjacent rooted edge. Strictly speaking,
we must use the black-white symmetry of the maps counted by $F_{\rm BMHP}$ which imposes
that there are in average as many black-empty vertices as white-empty ones.

A straightforward, though tedious exercise allows by a Lagrange inversion to get
a closed formula for the coefficient $G_{2n}$ of $g^{2n}$ in $G_{\rm BMHP}(g,z)$ 
(note that $G_{\rm BMHP}$ has only even powers of $g$ as bicubic maps always have
an equal number of black and white vertices). Details can be found in Ref.\HARM, with the
result
\eqn\exactenum{ G_{2n}(z)={3\over 2} {2^n\over (n+1)(n+2)} \sum_{0\leq 2p\leq j\leq n}
\left(-{1\over 2}\right)^p {2n-j \choose n}{n-j \choose p}{4n-2j\choose j-2p} z^j}

\subsec{Phase diagram}

A nice outcome of the above algebraic solution is the existence of a critical line $g_c(z)$ 
governing the large $n$ growth of $G_{2n}$, namely 
\eqn\Ggrowth{G_{2n}\sim {g_c(z)^{-2n}\over n^{5/2}}}
for generic values of $z$. The value of $g_c(z)$ is easily obtained from the singularity 
of $G_{\rm BMHP}$ with the result
\eqn\gc{g_c^2(z)=\left\{\matrix{\displaystyle{{1\over 8z}-{1\over 4 z^2}}& {\rm for}\ z\geq z_{+}\cr
\displaystyle{{1+8u+10u^2\over 8(1+2u)^8}}\ {\rm if}\ z=4u(1+2u)^4 & {\rm for}\ z_{-}
\leq z\leq z_{+}\cr} \right.}
Physically, the quantity ${\rm Log}\,g_c(z)$ may be interpreted as a free energy per vertex.
As explained in Ref.\CRI, the two values $z_{-}$ and $z_{+}$ correspond to two
special points 
\eqn\mcpoints{\eqalign{
g_{-}&={5^3 \sqrt{15} \over 2^{10}},\quad z_{-}=-{2^9\over 5^5}\cr
g_{+}&={\sqrt{15} \over 2^6},\quad z_{+}=2^5 \cr}}
at which the model becomes ``tricritical" with a large $n$ behavior 
\eqn\Ggrowthmc{G_{2n}\sim {g_{\pm}^{-2n}\over n^{7/3}}}
The first tricritical point ($g_{-}$) occurs at a (rather unphysical) negative value of $z$ and
corresponds to the so-called Lee-Yang edge singularity. The second tricritical point
($g_{+}$) occurs at a positive value of $z$ and corresponds to a change of determination 
of $g_c(z)$ in Eq.\gc, i.e.\ of the free energy. Physically, it corresponds to a 
{\it crystallization transition}
point above which particles tend to lie preferentially on vertices of a given color.
This may be seen by breaking the symmetry between black and white vertices upon introducing
different weights $z_\circ$ and $z_\bullet=z^2/z_\circ$ for particles occupying white or 
black vertices and by computing the difference of the densities of particles on the black and
the white sublattice. Letting $z_\circ\to z$ from the above yields a non zero density
difference for $z>z_+$ while the density difference vanishes for $z\leq z_+$. 
This transition lies in the universality class of the Ising model on random lattices.
This is corroborated by computing from Eq.\gc\ the order of the discontinuity of 
the free energy at $z_+$: one finds a discontinuity in the third derivative, immediately
translated into a ``thermal exponent" $\alpha=-1$, identified with that of the Ising transition 
on random graphs \ISING.

\newsec{Combinatorics of bicubic maps with hard particles via blossom trees: main results}
The aim of this paper is to provide an alternative derivation of the set of equations \matrecur\ above 
by a purely combinatorial approach based on the cutting of maps into so-called ``blossom trees". 
This approach is similar in spirit to the bijective methods developed in Ref.\SCHth\ using one-to-one
correspondences between, on the one hand, various classes of rooted planar maps and, on
the other hand, blossom trees, i.e.\ trees with appropriate ``charged" decorations (so-called buds 
and leaves) and simple charge characterizations. In the present case, a similar bijection also 
exists between rooted bicubic maps with hard particles and a particular class of blossom trees, 
which we shall refer to as {\it good} trees. These good trees however have no simple charge characterization 
and, as such, do not allow for a simple enumeration. Fortunately, we may rely on a more general 
many-to-one correspondence between what we shall call {\it admissible} trees, now with an easy 
charge characterization, and a larger class of bicubic maps with particles.  The generating function 
for the admissible trees with appropriate signed weights turns out to match exactly that of good 
trees, and therefore enumerates the desired rooted bicubic maps with hard particles.
This is the key for our combinatorial approach.

\subsec{Blossom trees}

The combinatorial enumeration of planar bicubic maps {\it without particles} involves a bijection with
so-called {\it blossom trees}. A blossom tree is a plane tree satisfying the following
properties:
\item{(1)} Its inner vertices are of two types: black or white, all of valence three.
\item{(2)} It is bicolored, i.e.\ all its inner edges connect only black inner vertices to white ones.
\item{(3)} It has two kinds of leaves: the so-called {\it buds}, connected only to black inner vertices 
and the so-called {\it leaves}, connected only to white inner vertices.
We attach a {\it charge} $-1$ to each bud and a charge $+1$ to each leaf.
\item{(4)} Its total charge ($\#$leaves $-$ $\#$buds) is equal to $3$.
\par
\noindent Note that this last property is equivalent to imposing that there is exactly one more
white vertex that black ones in the tree. Moreover, a direct consequence of the defining
properties (1)-(4) is that each inner edge of a blossom tree separates it into a piece 
of total charge $q_\circ\equiv 2{\rm \ mod\ }3$
starting with a white vertex, and a piece of total charge $q_\bullet=3-q_\circ\equiv 1{\rm \ mod\ }3$
starting with a black vertex. Such an edge will be called of type $(q_\circ:q_\bullet)$.
\fig{An example of blossom tree (a).  Here and throughout the
paper, buds are represented by black arrows and leaves by white arrows. 
The tree may be closed into a bicubic map (d) as follows:
Each pair of counterclockwise consecutive bud-leaf is glued
into an edge (b). The process is repeated (c) until one is left with exactly
three unmatched leaves. The latter are connected (d) to an additional
black vertex, thus producing a planar bicubic map. Choosing one of the three
unmatched leaves as a root for the blossom tree (crossed in (c)) amounts to
rooting the map (white arrow in (d)).}{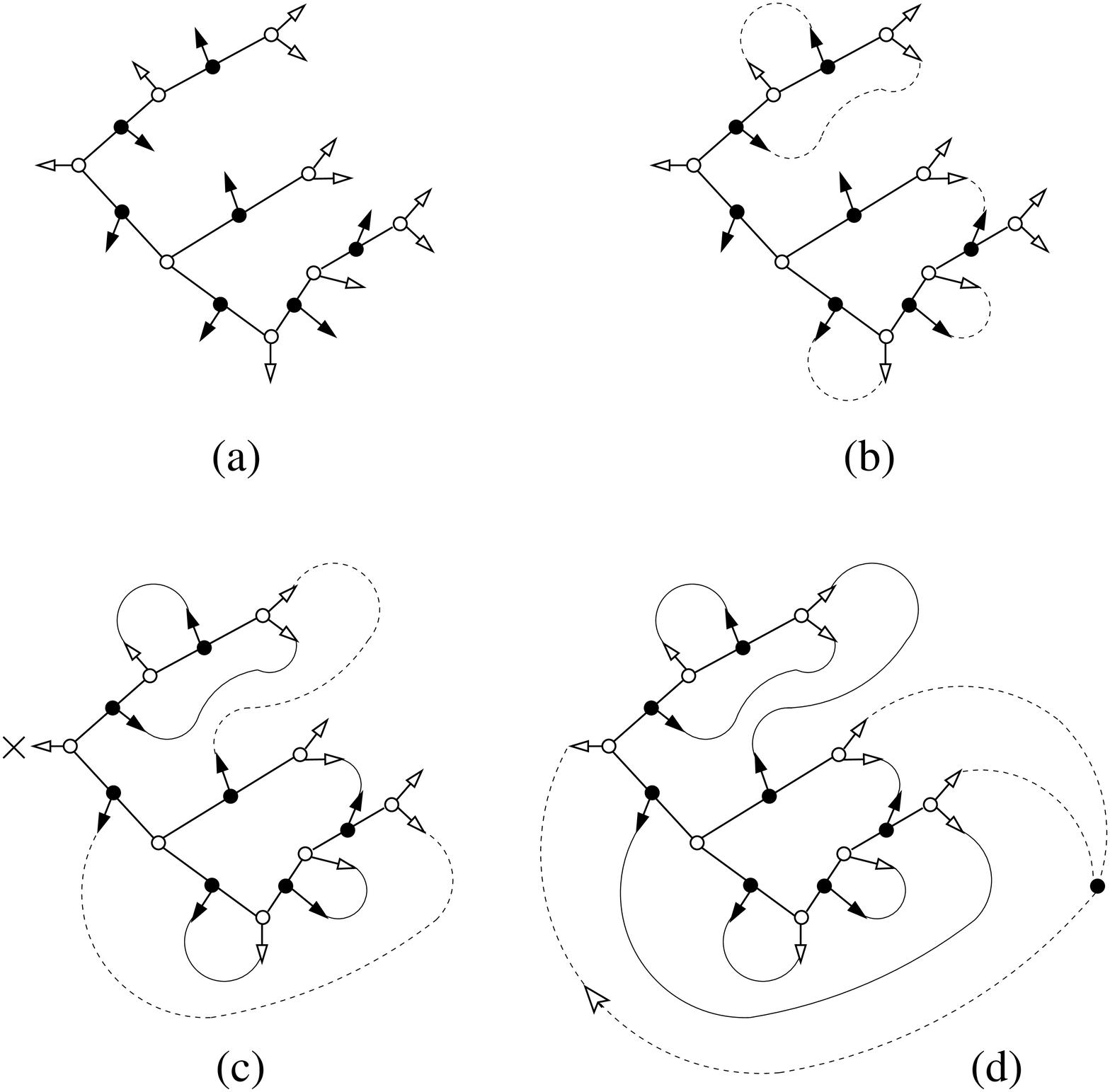}{11.cm}
\figlabel\refermeture
A blossom tree may be closed into a map as follows: we consider the counterclockwise sequence 
of buds and leaves around the tree and decide to connect all pairs made of a bud immediately 
followed by a leaf so as to form an edge (see Fig.\refermeture). Iterating this procedure, 
we end up with exactly three {\it unmatched} leaves, as the total charge of the tree is $3$. 
These unmatched leaves are then connected to an additional black vertex in the external face.
The net result is clearly a planar bicubic map. In the following, we will make use of {\it rooted}
blossom trees, i.e.\ trees with a distinguished leaf serving as root. More precisely, 
we shall always pick the root among the three unmatched leaves of the tree. Via the above
closing procedure, this allows to generate a rooted bicubic map, i.e. a map with a
marked oriented edge, namely that pointing from the additional black vertex to the
chosen unmatched leaf.

Conversely, any bicubic map may be cut into a blossom tree. The precise cutting procedure
will be recalled in Sect.4.1 below. The blossom trees thus obtained form a restricted
class, subject to specific restrictions which ensure 
that the charge is evenly distributed along the tree. More precisely, 
an edge of type $(q_\circ:q_\bullet)$ is said to be {\it regular} if and only if
it satisfies $q_\circ \geq 0$ and $q_\bullet \leq 1$. Any edge that is not regular is called 
{\it non-regular}. 
The trees obtained by cutting bicubic maps satisfy the additional constraint that all 
their inner edges are regular\foot{For the blossom trees at hand, this property simply amounts
to imposing that there is exactly one bud attached to each black vertex. The tree of
Fig.\refermeture-(a) satisfies this property.}. This furthermore leads to a bijection between rooted planar 
bicubic maps and blossom trees with only regular inner edges, rooted at one of their three unmatched
leaves [\xref\CONST,\xref\COL]. In this framework, the cutting/closing procedures are moreover 
inverse of one another.

\subsec{Admissible trees and main theorem}

We would like to extend the above bijection to the case of bicubic maps with hard particles.
A natural candidate for the trees to be considered are blossom trees satisfying (1)-(4) above,
with all their inner edges regular, and now carrying particles satisfying the hard particle constraint.
In the following, we shall denote by {\it HP} any edge that satisfies the hard particle constraint 
and {\it NHP} (non-hard-particle) any doubly-occupied edge. In other words, it seems natural
to introduce ``HP-regular" blossom trees, i.e. blossom trees with particles whose inner edges 
are both HP and regular. Clearly, the
standard cutting procedure of bicubic maps, performed regardless of particles, will produce
such trees. Unfortunately, not all the HP-regular blossom trees are created that way. Indeed
for an arbitrary HP-regular tree, the closing procedure may produce a map with NHP edges if
two occupied vertices get connected in the process. If so, the tree at hand does not arise
from the cutting of bicubic maps with particles.
We may still define {\it good} trees as the blossom trees whose inner edges are all HP and regular 
and whose closing does not create NHP edges.
The (rooted) good trees are in bijection with (rooted) planar bicubic maps with hard particles. 
However, these trees have no simple local characterization and cannot serve as such to enumerate the maps. 

This is the reason why we are led to consider a larger class of blossom trees, which we
shall call {\it admissible}.  As opposed to good trees, the hard particle condition 
need not be fulfilled on admissible trees and edges may be NHP on the tree. The key point
will be however to correlate the  NHP character of inner edges to their non-regularity.  
More precisely, we define an admissible tree as a blossom tree satisfying properties
(1)-(4), further completed by the condition that each inner vertex may be occupied or
not by a particle, and with the two extra conditions:
\item{(5)} All HP edges of the tree are regular.
\item{(6)} All NHP edges of the tree are non-regular. 
\par
\noindent Thanks to properties (1)-(4), these two conditions may be rephrased into
\item{(5')} Every HP edge is of type $(q_\circ:q_\bullet)$ with $q_\circ\geq 2$, 
hence $q_\bullet \leq 1$. 
\item{(6')} Every NHP edge is of type $(q_\circ:q_\bullet)$ with $q_\circ\leq -1$, 
hence $q_\bullet \geq 4$.
\par
\fig{An example of (a) admissible and (b) non-admissible blossom tree. 
As in Fig.\bicub, black-empty, white-empty, black-occupied and white-occupied inner vertices
are represented respectively by $\bullet$, $\circ$, $\circbullet$ and $\circcirc$. 
For both trees (a) and (b), the total charge is $3$ hence there are exactly three unmatched
leaves, here marked by a cross.  For each inner edge, we have indicated the value $q_\circ$ of the total 
charge for the part of the tree lying on the side of its white (empty or occupied) adjacent vertex
(the reader may immediately recover the value of $q_\bullet=3-q_\circ$). 
All edges satisfy $q_\circ\equiv 2 {\rm \ mod\ }3$. An edge is regular when $q_\circ \geq 2$.
In (a), all HP edges are regular and all NHP edges are non-regular, hence the tree is admissible. This
is not the case in (b) where the indicated edge is both regular and NHP.}{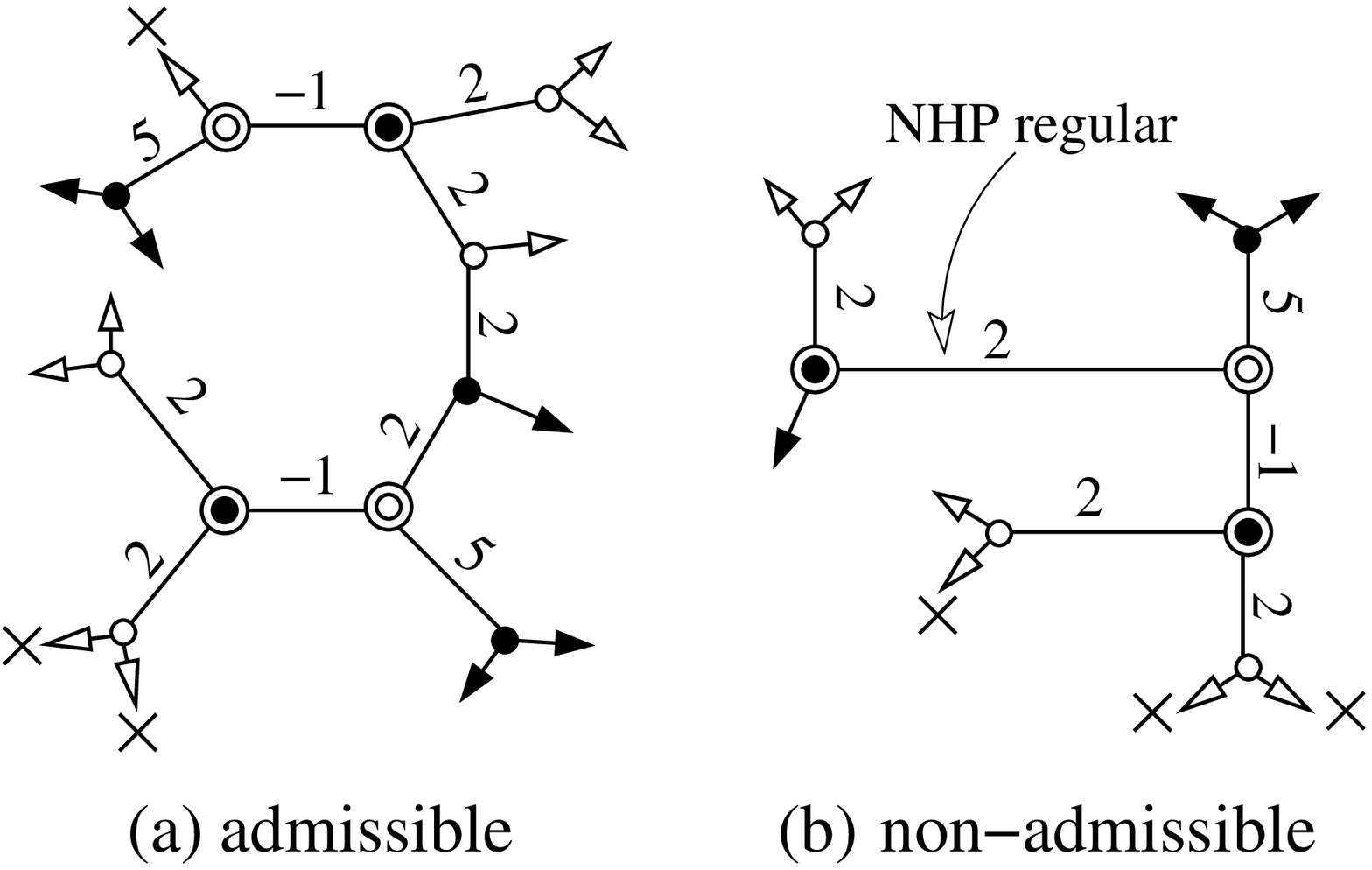}{10.cm}
\figlabel\admissible
\noindent Examples of admissible and non-admissible blossom trees are displayed in Fig.\admissible.
As before, we will make use of rooted admissible trees, i.e.\ admissible trees  
with a distinguished leaf among their three unmatched ones (marked by a cross on the examples of 
Fig.\admissible).

In the absence of NHP edges, the requirement (5) states that all edges are regular.
The good trees are therefore a particular subclass of the admissible trees. As it turns out,
the former may be obtained from the latter via an inclusion-exclusion principle which
allows to enumerate the good trees by counting admissible trees with appropriate signed
weights. This is our main result, stated as follows.
\medskip
\noindent {\underbar {Main theorem}}

\noindent For any admissible tree ${\cal T}$, we denote by $n({\cal T})$ the number of its NHP edges. 
We have the following sum rule: 
\eqn\mainth{\sum_{{\rm rooted\ admissible}\atop {\rm trees\ }{\cal T}}(-1)^{n({\cal T})}=
\sum_{{\rm rooted\ bicubic\ maps}\atop {\rm with\ particles\ }{\cal M}} 1 =
\sum_{{\rm rooted\ good}\atop {\rm trees\ } {\cal T}} 1} 
Here the first and last sums run over all (admissible, resp. good) trees with fixed numbers 
$n_\bullet$, $n_\circ$, $n_\smallcircbullet$ and $n_\smallcirccirc$ of black-empty, white-empty, 
black-occupied and white-occupied vertices. These trees are rooted at one of their three unmatched
leaves. The middle sum runs over planar maps with the same
numbers $n_\circ$, $n_\smallcircbullet$ and $n_\smallcirccirc$ of white-empty, black-occupied
and white-occupied vertices and with $n_\bullet+1$ black-empty vertices. These maps are rooted
maps, with a marked oriented edge originating from one of these $n_\bullet+1$ black-empty vertices. 

As the admissible trees have a simple charge characterization, their enumeration, with
the appropriate signed weight, is easily performed via a recursive construction.
The above sum rule allows to directly re-interpret this enumeration as counting 
rooted bicubic maps with hard particles.

\subsec{Local environments}

\fig{Inner edges adjacent to a white-empty vertex are HP, regular of type 
$(q_\circ=2:q_\bullet=1)$, and connect it to a black-empty or -occupied 
vertex.}{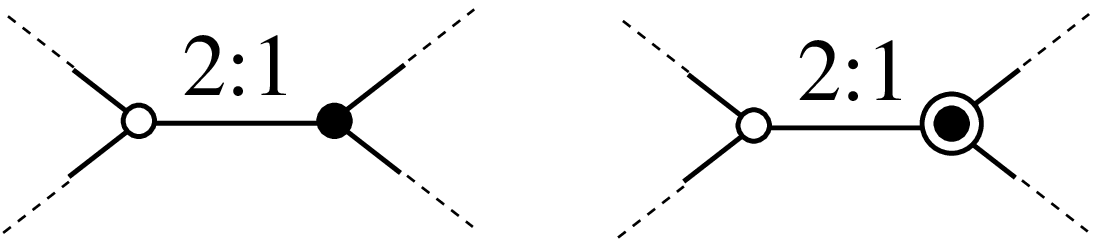}{8.cm}
\figlabel\envwhiteempty
In order to enumerate admissible trees, it proves useful to first inspect
all possible environments of an edge or a vertex on the tree and to characterize
them in terms of charges. The above characterizations (1)-(6) of admissible trees 
strongly limit these possible charge environments. A first limitation concerns 
the local environment of white-empty vertices. Indeed, any such vertex is adjacent to
either leaves, with charge $+1$, or inner HP edges, which are thus regular with 
$q_\bullet \leq 1$. As the vertex is trivalent and the total charge of the tree is 
$3$, we immediately deduce that only $q_\bullet=1$ is possible. Any inner edge attached 
to a white-empty vertex is therefore of type $(q_\circ=2:q_\bullet=1)$ (see Fig.\envwhiteempty). 

\fig{The four possible environments around a NHP edge. The NHP edge is necessarily of type
$(q_\circ=-1:q_\bullet=4)$. All other inner edges are HP regular of type
$(q_\circ=2:q_\bullet=1)$ except for the edge connecting $w$ to the black-empty vertex
with two buds, which is HP regular of type $(q_\circ=5:q_\bullet=-2)$.}{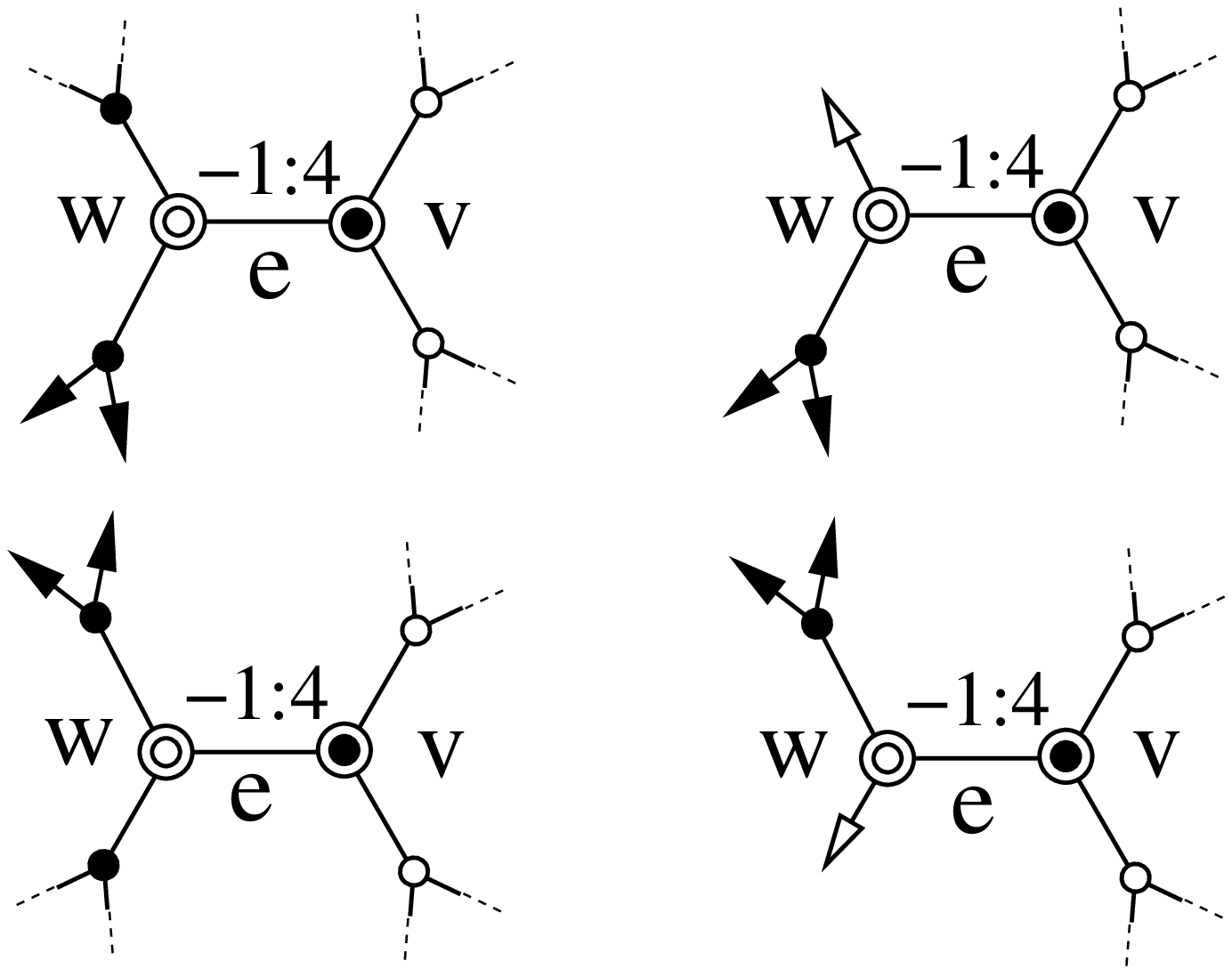}{10.cm}
\figlabel\envNHP
Let us now consider the environment of a NHP inner edge $e$ connecting a black-occupied 
vertex $v$ to a white-occupied vertex $w$ (see Fig.\envNHP). From (6'), the charge $q_\bullet$ 
of the piece of the tree attached to $v$ satisfies $q_\bullet \geq 4$. This piece is made of 
two parts, each attached to $v$ by one of its two remaining adjacent edges (other than $e$). 
As their respective charges add up to $q_\bullet \geq 4$, at least one of these two parts 
must carry a charge greater or equal to $2$, hence the corresponding edge adjacent to $v$ 
must be a regular inner edge. This edge connects $v$ to a white vertex, and, as it is regular, 
it must be HP, hence this white vertex is empty. We deduce from the above characterization 
around white-empty vertices that the charge of this part is exactly $2$.  The charge $q_\bullet -2$
of the other part is thus greater or equal to $2$ and from the same argument, it must be $2$ 
as well, starting from a regular inner edge connecting $v$ to a white-empty vertex. 
We finally deduce that $q_\bullet=4$, hence a complementary charge $q_\circ=-1$ for the piece 
of the tree starting from the white-occupied vertex $w$. This piece is also made of two parts, 
each attached to $w$ by one of the two remaining adjacent edges (other than $e$). One of 
these parts must carry a charge less than $-1/2$, hence it must be a regular edge connecting 
$w$ to a black-empty vertex. As all inner edges originating from a black-empty vertex are 
HP, thus regular with an attached charge larger than $2$, the only way to have a charge 
less than $-1/2$ is that this black-empty vertex be connected to two buds, with a total 
charge of $-2$. The second part attached to $w$ therefore has charge $1$, hence is either 
a leaf or a regular inner edge connecting $w$ to a black-empty vertex. To conclude, we find 
that all NHP edges are of type $(q_\circ=-1:q_\bullet=4)$ with one of the four environments 
displayed in Fig.\envNHP. 
\fig{Inner edges adjacent to a black-empty vertex are HP, regular of type 
 $(q_\circ=2:q_\bullet=1)$ or $(q_\circ=5:q_\bullet=-2)$. In this latter
case the black-empty vertex is connected to two buds and the white
vertex is occupied and adjacent to exactly one NHP edge}{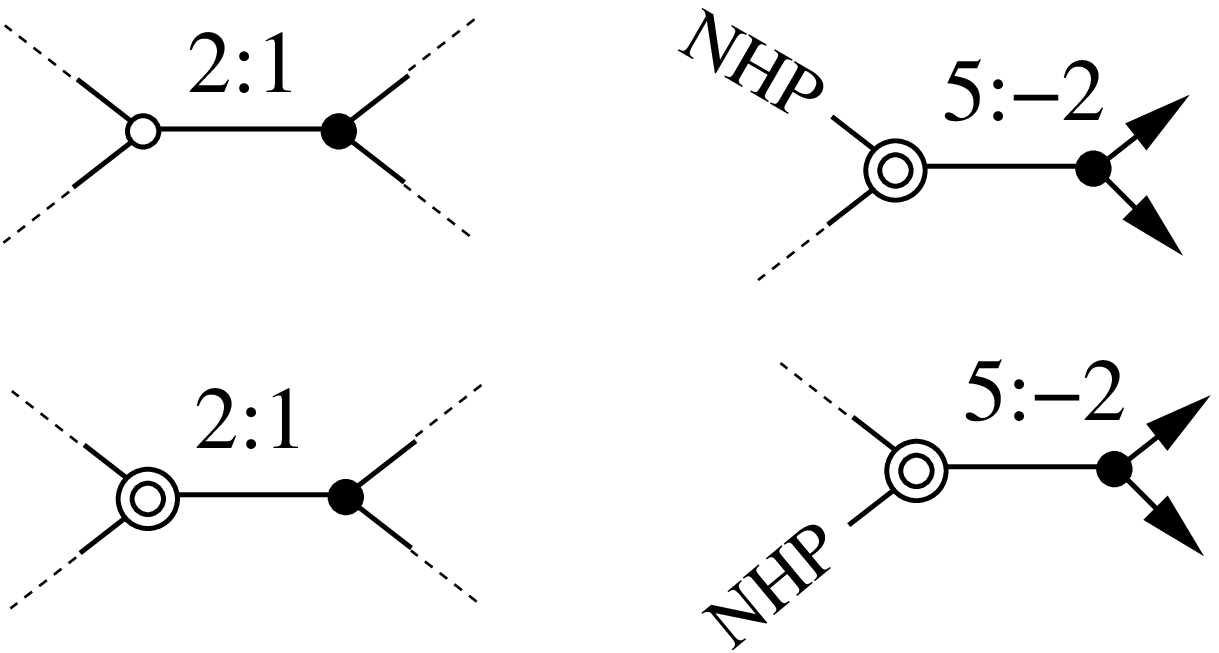}{8.cm}
\figlabel\envblackempty
If we finally consider HP, thus regular inner edges, we already know from the above 
analysis around white-empty vertices that these edges are of type $(q_\circ=2:q_\bullet=1)$ 
when the adjacent white vertex is empty, irrespective of the empty or occupied state
of the adjacent black vertex (see Fig.\envwhiteempty). It remains to consider the case 
of an inner edge connecting
a black-empty to a white-occupied vertex. Such an edge is clearly either of type 
$(q_\circ=2:q_\bullet=1)$, or of type $(q_\circ=5:q_\bullet=-2)$. In this latter case, 
the black-empty vertex is necessarily attached to two buds, while the white-occupied 
vertex is necessarily adjacent to exactly one NHP edge (see Fig.\envblackempty).
This exhausts all possible inner edge environments.

The vertex environments follow directly from these inner edge environments, plus the
possibility of attaching a leaf to white-empty or -occupied vertices and a bud
to black-empty or -occupied vertices. The vertex environments are moreover subject to
the constraint of total charge $3$ of the tree. 

\subsec{Recursive equations}

We are now ready for a combinatorial enumeration of admissible trees. More precisely,
we want to compute the generating function $G(g,z)$ for admissible trees rooted at
one of their three unmatched leaves, with weights
$g$ per vertex and $z$ per particle (i.e.\ a total weight of $gz$ per occupied vertex) 
and with {\it a weight $-1$ per NHP inner edge} so as to reproduce the generating function
for good trees according to Eq.\mainth. 

This enumeration is simplified by considering {\it planted} trees, i.e trees with a 
distinguished leaf or bud. As all leaves except for the three unmatched ones are paired
with buds, we may write: 
\eqn\gleafbud{G(g,z)= G_{\rm leaf}(g,z)-G_{\rm bud}(g,z)}
in terms of the generating functions $G_{\rm leaf}(g,z)$ (resp. $G_{\rm bud}(g,z)$) of 
admissible trees with a distinguished leaf (resp. bud). 
From the main theorem \mainth, this equation translates into
\eqn\gener{G_{\rm BMHP}(g,z)=g(G_{\rm leaf}(g,z)-G_{\rm bud}(g,z))}
where $G_{\rm BMHP}(g,z)$ is as before the generating function
for planar bicubic maps with hard particles, with a marked oriented edge originating
from a black-empty vertex, and with weights $g$ per vertex and $z$ per particle.

The characterization of admissible trees involves charges $q_\bullet$ and $q_\circ$ for 
pieces of the tree on both sides of a given inner edge. We shall refer to such pieces as 
{\it half-trees}. For convenience, we may also consider buds and leaves as half-trees, 
in which case we may decide to see the half-trees as dangling from the inner vertex 
to which they are attached. For instance, any half-tree attached to a white-empty vertex 
is of charge $+1$, whether it reduces to a leaf or it starts by a regular inner
edge (see Fig.\envwhiteempty). We may then introduce the following generating functions: 
\item{(i)} $V$ for half-trees of charge $1$ attached to a white-empty vertex
\item{(ii)} $S$ for half-trees of charge $1$ attached to a white-occupied vertex
\item{(iii)} $T$ for half-trees of charge $4$ attached to a white-occupied vertex
\item{(iv)} $R$ for half-trees of charge $2$ attached to a black-empty vertex
\par
\fig{A pictorial representation of the recursive relations (3.4) for half-trees, 
respectively attached to a white-occupied vertex, with charge $1$ ($S$) or $4$ ($T$),
to a black-empty vertex, with charge $2$ ($R$) and to a white-empty vertex,
with charge $1$ ($V$). The various terms collect all possible environments
compatible with the charge characterization of admissible trees. We also
indicate the various weights $g$, $z$ and minus signs.}{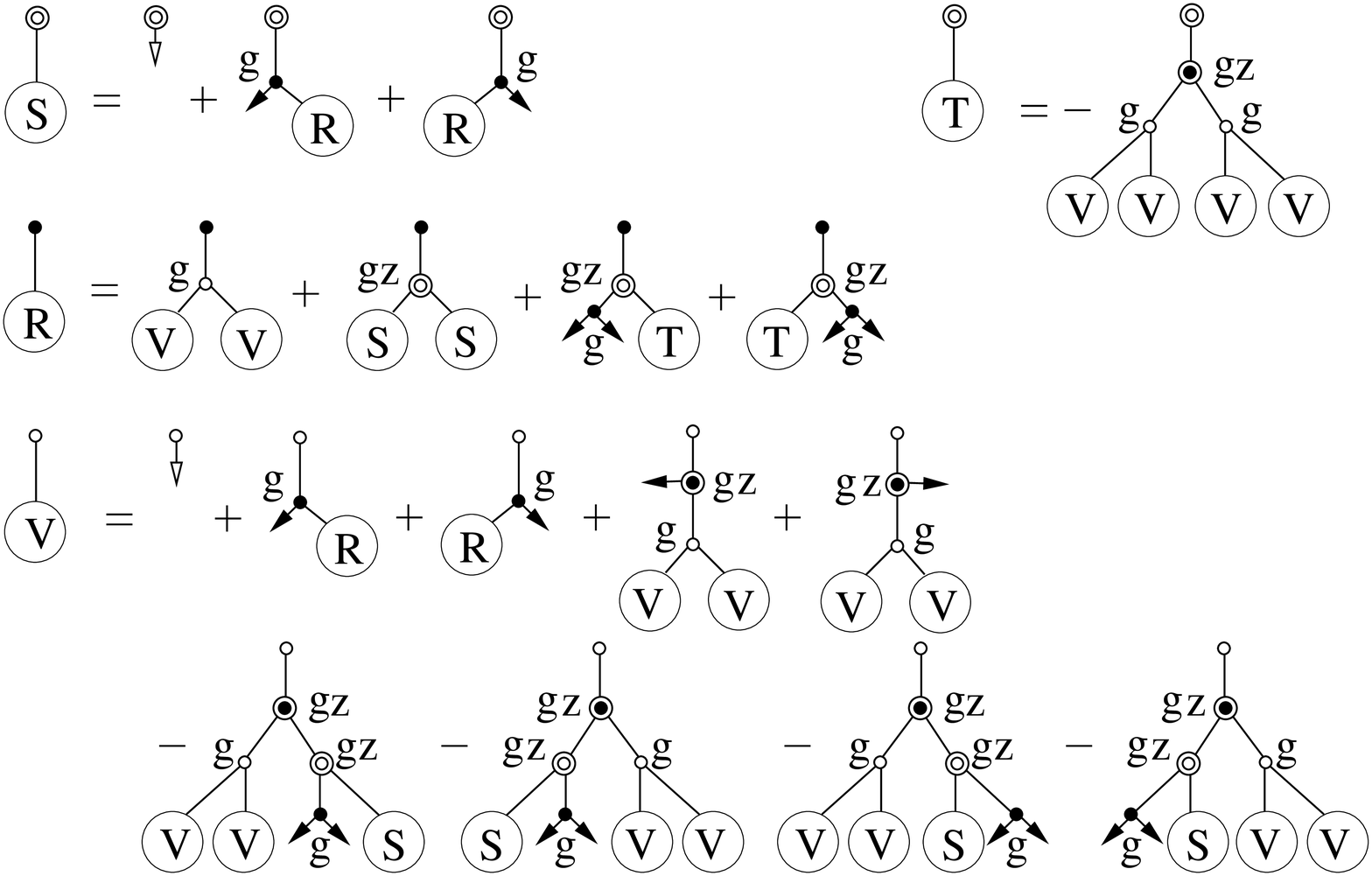}{13.cm}
\figlabel\figrecur
\noindent By inspecting the environments of the root vertex of these half-trees, 
we find the following recursion formulas, easily read off Fig.\figrecur:
\eqn\recur{\eqalign{
S&= 1 + 2 g\, R \cr T &= -g^3z\, V^4 \cr R &= g\, V^2+gz\, S^2 +2 g^2z\, T\cr
V&= 1+2 g\, R + 2 g^2z\, V^2 -4 g^4z^2\, V^2 S\cr}}
The last equation may be rewritten as $(1-2g^2z\, V)(V-S-2g^2z\, VS)=0$, hence it
is equivalent to 
\eqn\simpV{V=S+2g^2z\, VS}
These equations, here derived in a purely combinatorial way, match exactly those 
\matrecur\ obtained from the matrix model solution.
\fig{A pictorial representation of the expressions (3.6) for $G_{\rm leaf}$
and $G_{\rm bud}$, as obtained by collecting all possible environments
of the leaf or bud compatible with the charge characterization of admissible trees. 
We also indicate the various weights $g$ and $z$.}{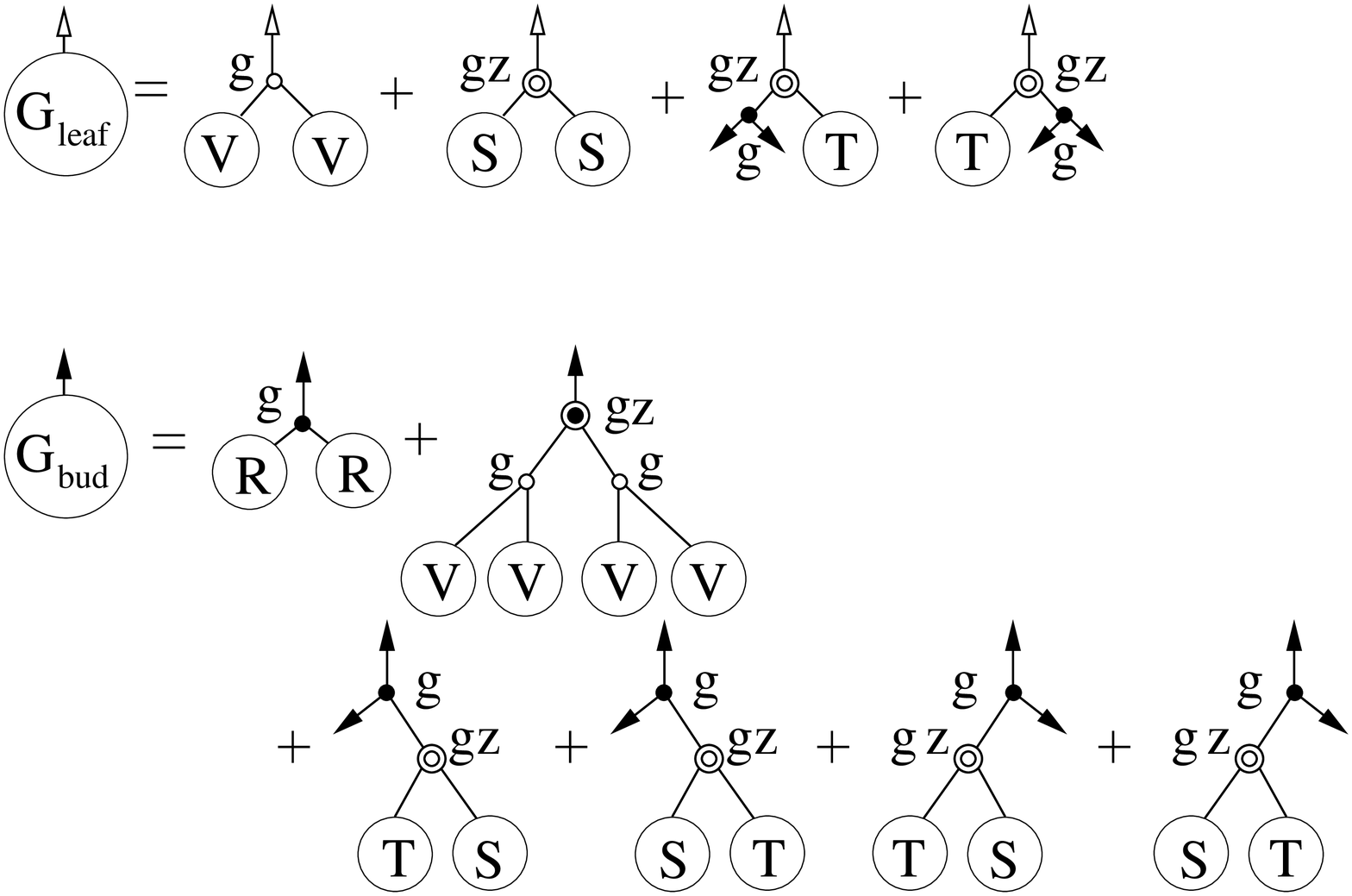}{11.cm}
\figlabel\planted
\noindent Similarly, by inspection of the environment of the distinguished leaf or bud in planted
trees, we immediately get the relations (see Fig.\planted):
\eqn\gVR{\eqalign{G_{\rm leaf}&= g\, V^2+gz\, S^2 +2 g^2z\, T = R\cr
G_{\rm bud}&= g\, R^2 +g^2\,z\,(g\,V^4+4\,S\,T)\cr}}
Equations \gener, \gVR\ and \recur-\simpV\ provide a closed system of algebraic equations
for the enumeration of planar rooted bicubic maps with hard particles. The first terms
of an expansion in $g$ of $G_{\rm BMHP}(g,z)$, counting rooted bicubic maps with, say, up 
to $10$ vertices, read for instance:
\eqn\firstterms{\eqalign{G_{\rm BMHP}(g,z)&= g^2 (1+z) +g^4 (3 +9z +3 z^2) 
+g^6 (12 + 60 z + 66 z^2 + 12 z^3) \cr & +g^8 (56 + 392 z + 780 z^2 + 460 z^3 + 56 z^4)\cr 
& +g^{10} (288 + 2592 z + 7584 z^2 + 8400 z^3 + 3168 z^4 + 288 z^5)+{\cal O}(g^{12})\cr}}
These values match exactly the general formula \exactenum. The free energy $F_{\rm BMHP}$
is easily recovered from $G_{\rm BMHP}$ by inverting Eq. \FtoG\ into
\eqn\GtoF{F_{\rm BMHP}(g,z)={2\over 3}\int_1^{\infty} {du\over u} G_{\rm BMHP}\left(
{g\over u},z u\right)}

\newsec{Cutting of rooted bicubic maps with non-sectile edges: generalities}

The main theorem \mainth\ is a consequence of a general relation between maps
and trees obtained via some cutting procedure. In this section, we specialize
to the case of planar bicubic maps without particles, but with possibly a number
of distinguished {\it special} edges, which we prevent from being cut in the process.
This general framework will be applied in the next section to the particular
case of bicubic maps with hard particles by relating the special edges to NHP
ones.

For starters, we first recall the cutting procedure in the standard case 
of bicubic maps with no special edges [\xref\CONST,\xref\COL]. We then extend the analysis in the
presence of special {\it non-sectile} edges.

\fig{Cutting procedure for a typical rooted bicubic map. The various steps (a)-(d)
are detailed in the text.}{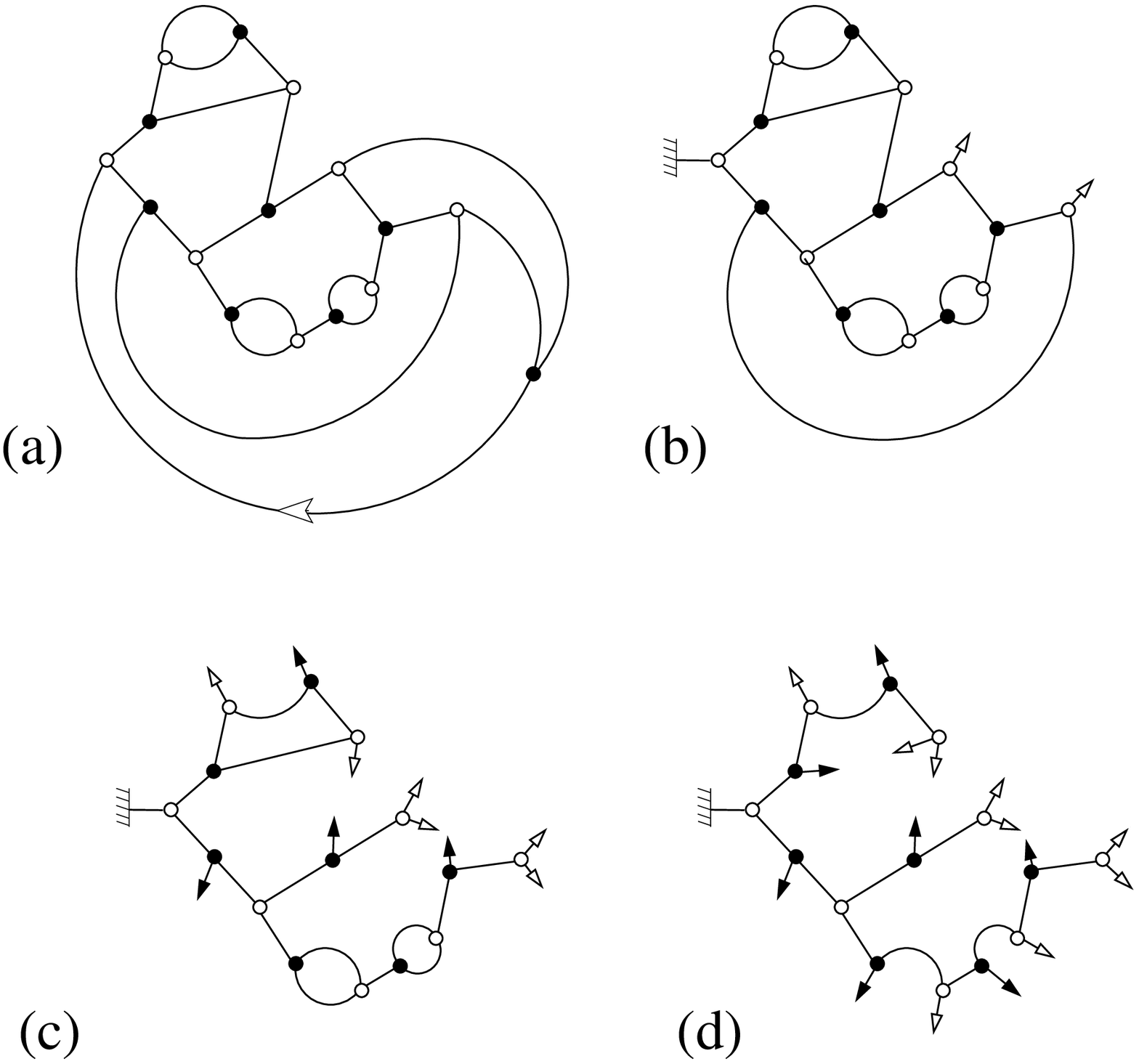}{9.cm}
\figlabel\decoupage

\subsec{Preamble: Cutting of ordinary bicubic maps}

The cutting of a planar rooted bicubic map into a blossom tree was discussed in Refs.[\xref\CONST,\xref\COL]. 
As illustrated in Fig.\decoupage, it is performed as follows:
\item{(a)} We draw the map in the plane with the external face lying to the left
of its root edge (marked oriented edge originating from a black-empty root vertex) 
(see Fig.\decoupage-(a)). 
\item{(b)} We erase the root vertex and replace the three adjacent edges by
leaves. The leaf corresponding to the rooted edge is marked and taken as a root
(represented as a grounding pictogram in Fig.\decoupage-(b)).
\item{(c)} Starting from this root, we visit successively each edge adjacent to the new
external face in counterclockwise direction. The edge is cut into two half-edges if 
(i) it originates from a black vertex and (ii) the cutting does not disconnect the 
remaining graph. The half-edge connected to the adjacent black vertex is replaced by a bud,
while that connected to the adjacent white vertex by a leaf. This results in merging 
a number of faces of the original map with the external one.
\item{(d)} We iterate the cutting procedure, each time by visiting counterclockwise 
the edges adjacent to the newly produced external face. We stop the process once
all faces have been merged and the remaining graph is a tree (see Fig.\decoupage-(d)).
\par
\noindent This cutting procedure is an {\it injection} from maps to blossom trees. Indeed, we
may easily re-build the original map from its associated tree by the standard closing procedure
discussed in Sect.3.1 (see Fig.\refermeture) [\xref\CONST,\xref\COL].
The injection is from maps with a rooted edge (originating from a black vertex) 
to blossom trees with a root taken among their three unmatched leaves. 

\fig{Equivalent formulation of the cutting procedure using geodesic distances. In
(a), we label each face by its geodesic distance from the external face (labeled $0$).
Representing the seam  by a wavy line, the geodesic tree, made of all leftmost minimal 
paths is drawn in dashed lines. The edges dual to those of this tree are replaced
in (b) by bud-leaf pairs. The blossom tree is finally obtained by erasing 
the root edge and root vertex (c), replacing them by the root of the 
tree.}{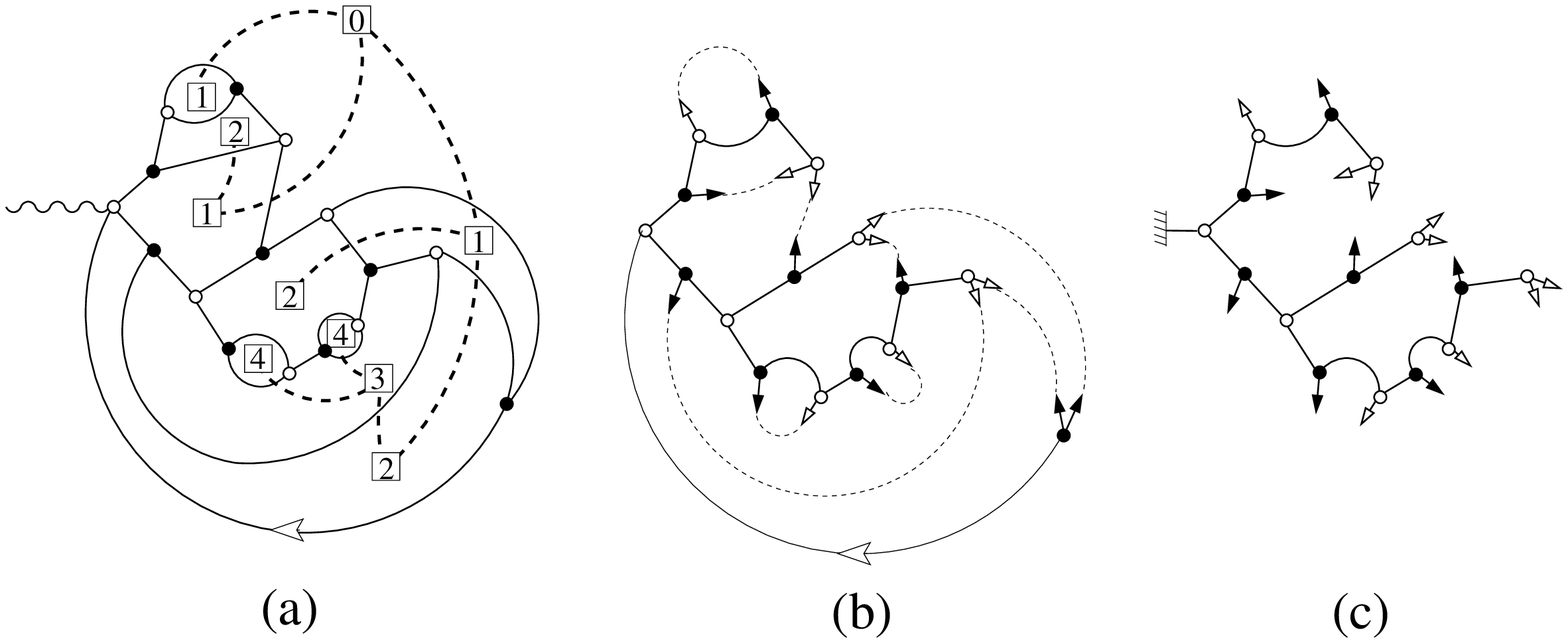}{13.cm}
\figlabel\decoupagetwo
An important property of the above cutting is that it may be reformulated
by use of {\it geodesic distance} between faces of the map \SCHth. Taking the
external face as origin, we may define its geodesic distance to any
other face of the map as the minimal number of edges to be crossed by a path 
stepping from face to face from the external face to the face at hand.
The paths are further required (i) to leave black vertices on their left,
and (ii) not to cross a seam connecting the extremity of the
rooted edge to infinity. Such paths realizing the minimal 
number of steps are called minimal paths. Any two such paths $P_1$ and $P_2$  
may be ordered by comparing their first differing steps $s_1$ and $s_2$. The
path $P_1$ is said to lie on the left of $P_2$ if $s_1$ lies on the left of $s_2$,
looking toward the map from infinity. This ordering allows to assign to each
face its {\it leftmost minimal path}. The union of these paths for all faces
forms a ``geodesic tree" (see Fig.\decoupagetwo-(a)). The above cutting procedure 
may be rephrased by cutting the edges dual to those of the geodesic tree and 
replacing them by bud-leaf pairs (see Fig.\decoupagetwo-(b)). The original root
edge and root vertex are removed together with the attached
buds and replaced by a distinguished leaf (see Fig.\decoupagetwo-(c)).  

The image of bicubic maps under cutting has the simple charge characterization 
that all edges be regular, as will be re-derived below. This allows for a nice bijective 
census of these maps.

\subsec{Cutting procedure in the presence of non-sectile edges}

\fig{A typical example of cutting of a rooted bicubic map (a) with special,
non-sectile edges, represented by thickened grey lines. In (b), we have erased the root
vertex of the map and replaced the adjacent edges by leaves, the extremity of the
rooted edge serving as a root (grounding pictogram). We iterate the standard 
cutting procedure (c)-(e) without cutting the special edges until a tree is
obtained (e). That the cutting ends up with a tree is guaranteed by the absence
of counterclockwise loops of either special or white-to-black edges. All non-special
edges are regular but special edges may be regular or non-regular (f).}{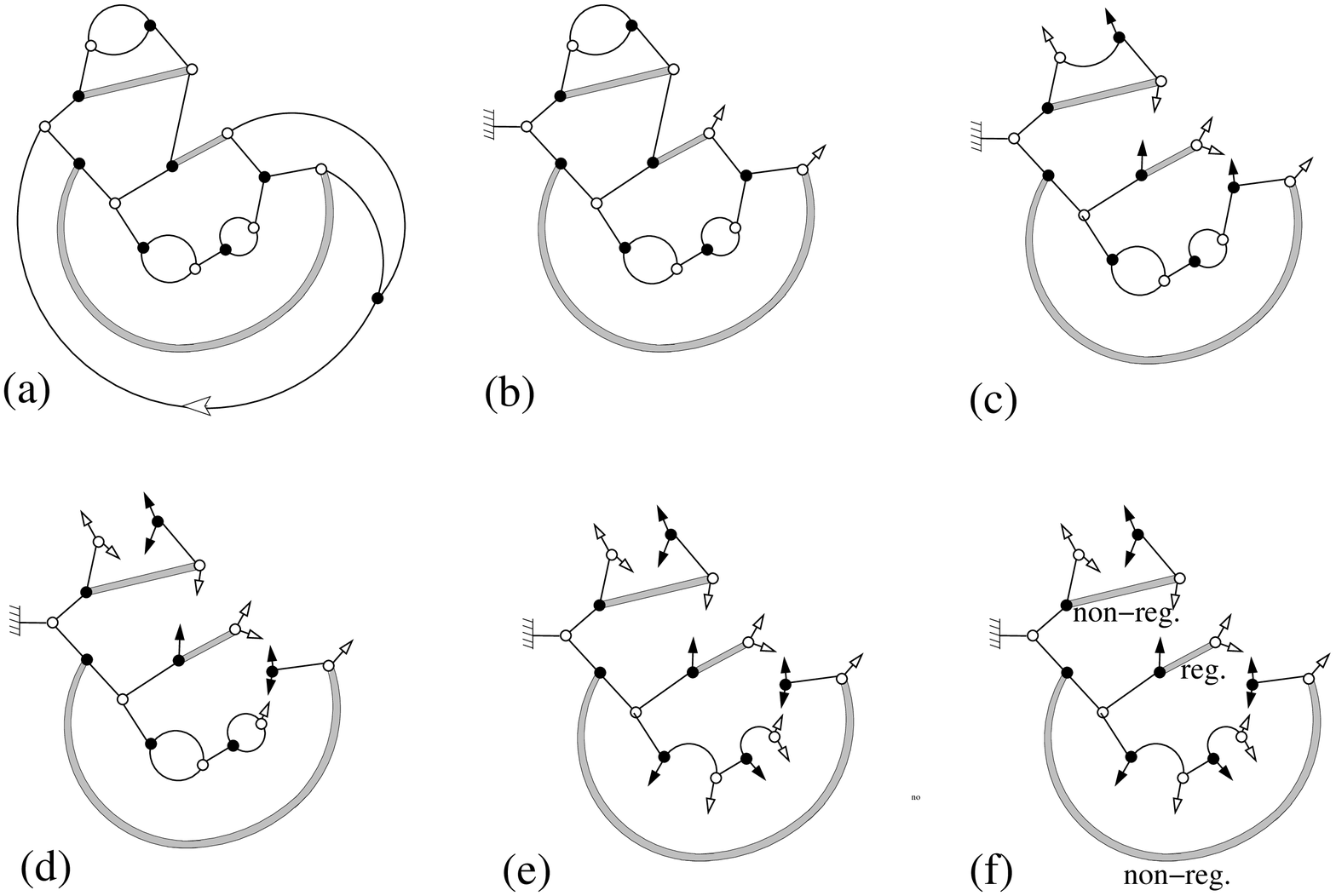}{13.cm}
\figlabel\decmarquage
The above cutting procedure is a particular case of a more general cutting
which also allows for the presence of {\it non-sectile} edges within the map. More
precisely, we may pick in a bicubic map a number of {\it special} edges 
which we decide {\it a priori} not to cut. Edges which are not special will be called
{\it normal}. We can repeat the above cutting process
following the same rules (a)-(d) of Sect.4.1, but with the extra requirement that special 
edges should not be cut. There is no guaranty in general that the cutting process may be carried on 
so as to produce a tree. In the following, we shall therefore restrict ourselves to 
{\it acceptable} choices of special edges, i.e.\ such that the cutting procedure 
does indeed produce a tree. In other words,
we forbid any choice of special edges for which a face of the map would be surrounded
only by special edges and white-to-black edges counterclockwise. We also forbid special
edges to be adjacent to the root vertex of the map. From now on, 
we shall always consider maps with an acceptable choice of special edges. 

The presence of special edges allows for the definition of 
a {\it new geodesic distance} using paths with the same requirements as before,
namely that (i) these paths leave black vertices on their left and (ii) do not cross 
a seam connecting the extremity of the rooted edge to infinity and with the
additional requirement that (iii) they now avoid the special edges.
A choice of special edges is acceptable if and only if each face remains at
a finite distance from the external face. Moreover, the cutting procedure in
the presence of special edges may still be rephrased as cutting the edges
dual to those of the geodesic tree made of all leftmost minimal paths (from the
external face to all other faces) that avoid the special edges.
A typical example of cutting procedure with special edges is displayed in 
Fig.\decmarquage.

We now wish to characterize the trees obtained by the above cutting.
Clearly, by construction, any such tree satisfies the following properties:
\item{(S1)} Its inner vertices are of two types: black or white vertices, 
all with valence three. Its inner edges are of two types, normal or special.
\item{(S2)} It is bicolored, i.e.\ all its inner edges connect only black to white vertices.
\item{(S3)} It has two kinds of leaves: buds, connected only to black inner vertices 
and leaves, connected only to white ones.
\item{(S4)} Its total charge ($\#$leaves $-$ $\#$buds) is equal to $3$.
\par
\noindent These are nothing but properties (1)-(4), i.e.\ the tree is indeed a blossom
tree as defined in Sect.3.1. The tree further enjoys the following charge characterization:
\item{(S5)} Normal (i.e.\ non special) edges are regular.
\par
\noindent From properties (S1)-(S4), this last property may be rephrased into:
\item{(S5')} Normal edges are of type $(q_\circ:q_\bullet)$ with $q_\circ\geq 2$ 
and $q_\bullet \leq 1$. 
\par
\fig{Relation between the charges $q_\bullet$ and $q_\circ$ of the two pieces
of tree ${\cal T}_1$ and ${\cal T}_2$ on both sides of an edge $e$, and the geodesic 
distances $m$ and $k$ from the external face of the two faces $\mu$ and 
$\kappa$ adjacent to $e$ on the map. We use the convention that the face 
$\mu$ follows the face $\kappa$ counterclockwise around the black endpoint of $e$.
We picked a particular configuration of the three unmatched leaves of the tree,
here all lying in ${\cal T}_1$. We read from the picture 
$q_\circ=(k-p)-(m-p)=k-m$, hence $q_\bullet=3+m-k$.}{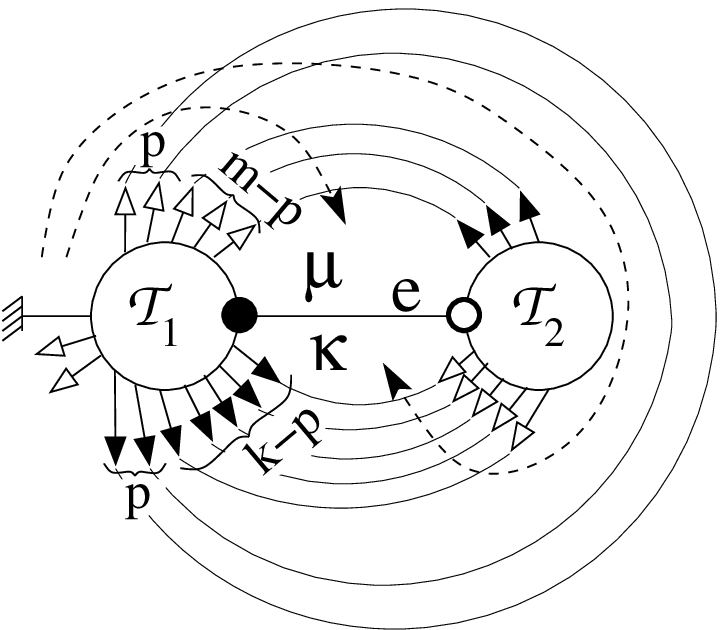}{7.cm}
\figlabel\thsample
\fig{Relation between the charges $q_\bullet$ and $q_\circ$ of the two pieces
of tree on both sides of an edge, and the geodesic distances $m$ and $k$ 
from the external face of the two faces adjacent to this edge on the map. We use the 
convention that the face at distance $m$ follows that at distance 
$k$ counterclockwise around the black endpoint of the edge. 
There are $12$ possible relative positions for the three unmatched leaves of the tree. 
In all cases on the left, we find $q_\circ=k-m$ while on the right $q_\bullet=m-k$.
The absence of leftmost minimal path crossing the edge imposes inequalities
on the geodesic distances $k$ and $m$ as shown. These 
inequalities translate into $q_\circ \geq 0$ or $q_\bullet \leq 1$.}{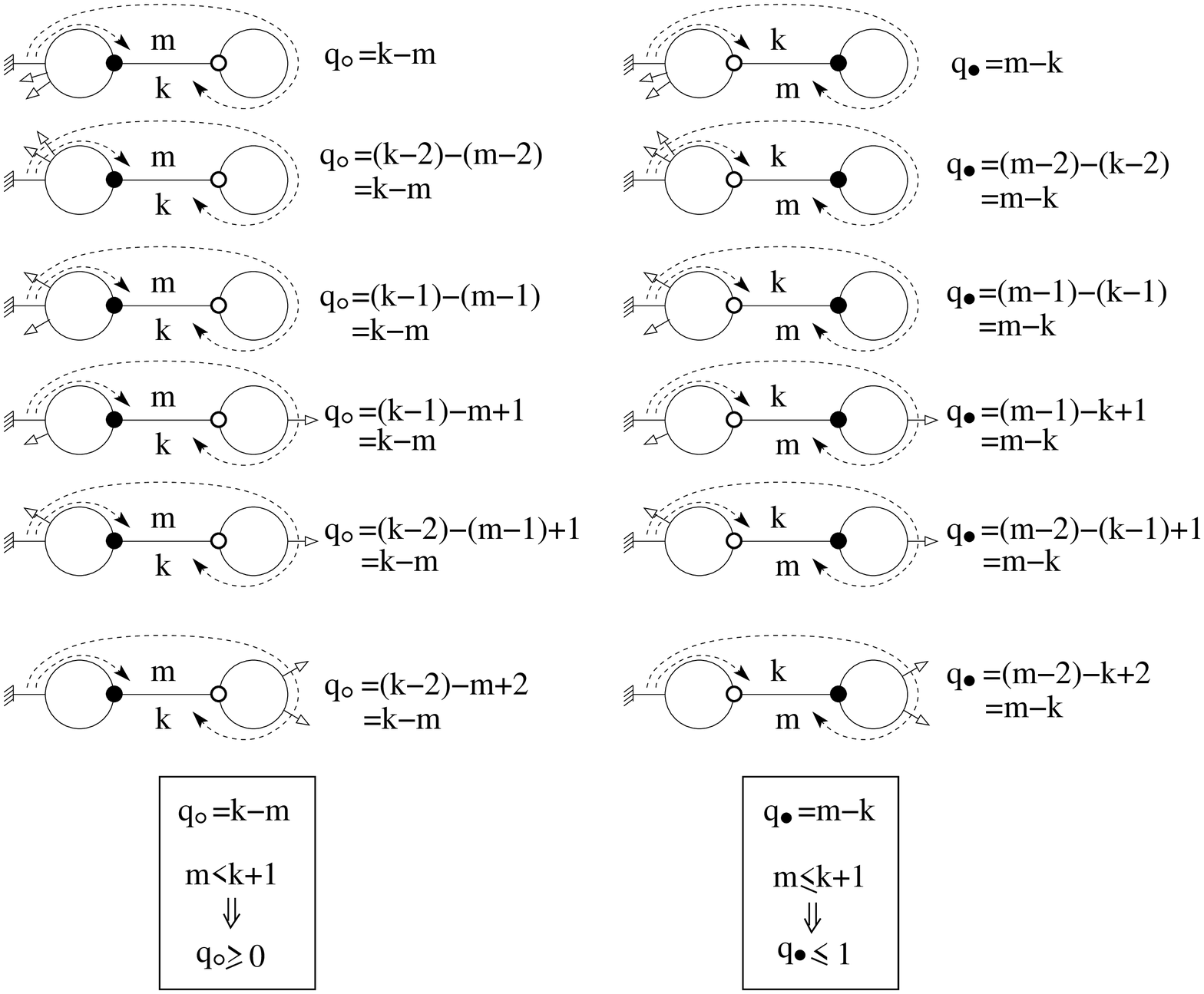}{13.cm}
\figlabel\thzero
\noindent The properties (S1)-(S4) are obvious. Let us now show the property (S5). Recall
that, as is true for any blossom tree, we have $q_\circ\equiv 2{\rm \ mod\ }3$ and 
$q_\bullet\equiv 1{\rm \ mod\ }3$
for the charge on both sides of any inner edge $e$. Moreover, these charges may
be related to geodesic distances as follows.
Let us denote respectively by $k$ and $m$ the geodesic distances of the two faces (resp. $\kappa$ and
$\mu$) adjacent to $e$ on the original map, with the convention
that the face $\mu$ follows the face $\kappa$  counterclockwise around the black
endpoint of $e$ (see Fig.\thsample). As a consequence of the leftmost minimal
path formulation of the cutting, we may relate the charges $q_\circ$ and $q_\bullet$
of the two pieces of the tree on both sides of $e$ to the geodesic
distances $k$ and $m$. Suppose for instance a situation like that of Fig.\thsample\ in
which the three unmatched leaves of the tree all lie on the same side of $e$, that
of its black endpoint. Denoting the piece of the tree on this side by ${\cal T}_1$ and
by ${\cal T}_2$ the piece on the side of the white endpoint of $e$, suppose moreover 
that the root is the last of the three unmatched leaves encountered by turning clockwise 
around ${\cal T}_1$ from $e$. Then the leftmost minimal path leading to $\mu$ will cut 
$m$ edges, resulting in the creation of $m$ leaves in ${\cal T}_1$ among which $m-p$ are 
matched to buds in ${\cal T}_2$, and $p$ are matched to buds of ${\cal T}_1$ by edges turning
around ${\cal T}_2$ (see Fig.\thsample). The leftmost minimal path leading to $\kappa$
has to cross these $p$ edges plus $k-p$ other edges, thus resulting in the
creation of $k-p$ leaves in ${\cal T}_2$ matched to $k-p$ buds in ${\cal T}_1$. All other
bud-leaf pairs correspond to cut edges not crossed by the two leftmost minimal
paths we are considering and thus lie on the same side of $e$ on the tree, giving
a net contribution zero to the charges of either piece. The charge of ${\cal T}_1$ is
therefore $q_\bullet=3+m-k$ while the charge of ${\cal T}_2$ is 
$q_\circ =k-m$.  There are $12$ possible relative positions of the
three unmatched leaves in either piece, as displayed in Fig.\thzero. A detailed 
case by case inspection of these $12$ cases shows that $q_\circ=k-m$ if the root 
of the tree lies on the side of the black endpoint of $e$ and $q_\circ=3+k-m$
if it lies on the side of the white endpoint. 

Finally, to prove (S5), we further assume that the edge $e$ is normal. The fact that it was 
not cut in the cutting process implies some inequality on the geodesic distances $k$ and $m$.
Indeed, we must have $m\leq k+1$ to eliminate the possibility of a shorter path 
leading to $\mu$ by crossing $e$ from $\kappa$. 
As the leftmost minimal path leading to $\mu$ lies on the right 
of that leading to $\kappa$ when the root lies on the side of the black endpoint of $e$, 
we must further forbid $m=k+1$ in this case to eliminate the possibility of
a path of the same length, but more to the left, by crossing $e$ (see Fig.\thzero).
These inequalities, together with the relations above between $k$, $m$ and
the charges $q_\circ$ and $q_\bullet$ immediately lead to the inequalities
$q_\circ \geq 0$ or $q_\bullet \leq 1$.  Combining this with $q_\circ=3-q_\bullet\equiv 2{\rm \ mod\ }3$, 
we deduce (S5'), or equivalently (S5), i.e $e$ is regular.
Note that there is no such characterization for special edges which may be regular 
or not (see Fig.\decmarquage-(f)).

To conclude this section, let us insist on the property (S5)(S5') which we rephrase as:

\noindent{\underbar {Proposition (C1):}} 
\noindent{\sl If a tree ${\cal T}$ is obtained by the cutting of  
a map ${\cal M}$ with a number of special (non-sectile) edges, then all normal 
({\rm a priori sectile}) edges of ${\cal T}$ are regular. Equivalently, if
an edge of ${\cal T}$ is non-regular, it must be a special edge.}

\fig{The cutting of a normal, thus necessarily regular edge, of ${\cal T}$ 
by a leftmost minimal path crossing that edge implies inequalities between 
the distance $k$ of the face adjacent to that edge (and such that the black endpoint
lies on the left when crossing the edge from that face) and the depth 
$m$ of the other adjacent face in the closing process of ${\cal T}$ 
(number of arches encircling that face). If we have picked an edge with minimal 
$k$, $k$ is also the depth of the first face in the closing process of ${\cal T}$
and we may relate $m$ and $k$ to the charges $q_\circ$ and $q_\bullet$ of the two 
halves of the tree as shown.  The inequalities then translate into the fact that 
the edge is not regular, a contradiction.}{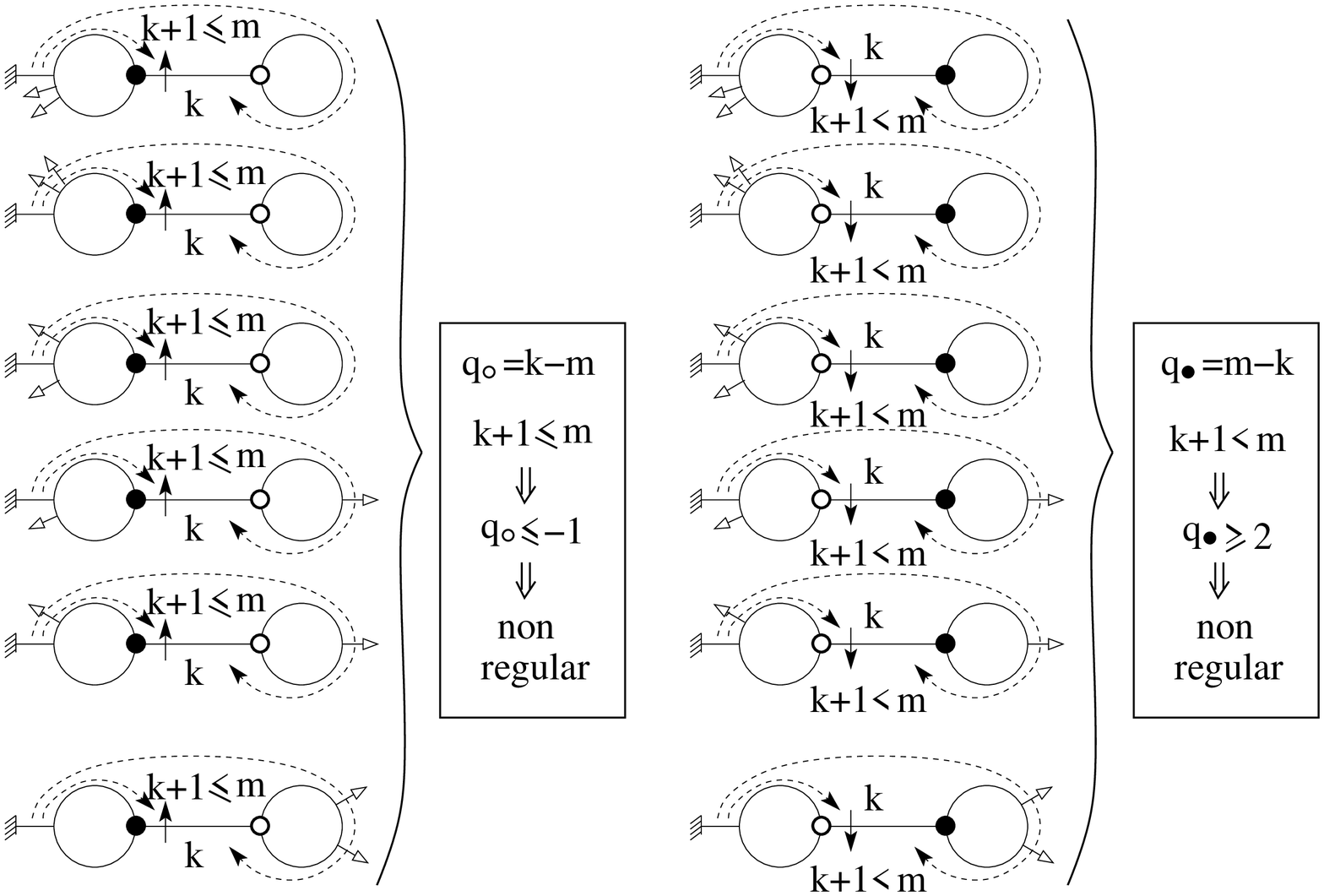}{13.cm}
\figlabel\thtwo
\subsec{Closing-cutting procedure: bijection}

The above cutting procedure sends any planar rooted bicubic map with
an acceptable choice of special edges onto a rooted blossom tree with special edges
and with the properties (S1)-(S5) above. This correspondence is clearly injective
as the map is easily restored from the tree by the same closing algorithm 
as that discussed in Sect.3.1. Moreover,
as we shall now see, the properties (S1)-(S5) entirely characterize the image
of the maps under cutting. We have therefore a bijection between, on one hand,
planar rooted bicubic maps with acceptable choices of special edges and, on
the other hand, blossom trees with special edges satisfying properties (S1)-(S5),
and with a distinguished leaf among their three unmatched ones.
That the closing of such a tree builds a rooted bicubic map with special edges is obvious. 
That the corresponding choice of special edges is acceptable is also clear as
each face of the map is attainable from the external face by a path whose edges
are dual to the closing arches plus possibly the duals of some edges adjacent to the root
vertex of the map. We are now left with the task of showing that closing
a tree ${\cal T}$  with the above properties into a map ${\cal M}$ with the
same special edges and re-opening it via the cutting procedure that preserves 
these special edges rebuilds ${\cal T}$ . It is clear that the cutting procedure 
produces a tree ${\cal T}'$ satisfying (S1)-(S5) and with the same number of 
inner vertices as ${\cal T}$, thus the same number of inner edges. To prove that ${\cal T}={\cal T}'$,
it is thus sufficient to prove that all the inner edges of ${\cal T}$ are in ${\cal T}'$.
It is clear that the special edges of ${\cal T}$ and ${\cal T}'$ are the same 
as they are not affected by the closing-cutting process. It is therefore sufficient 
to show that no normal edge of ${\cal T}$ may be cut in the cutting process. Assuming
the contrary, let us pick a normal edge $e$ of ${\cal T}$, supposedly 
cut in the cutting process and with {\it minimal distance} in ${\cal M}$ from the 
external face. By such edge distance, we mean the distance of that of the two faces,
say $\kappa$, adjacent to $e$ such that the black vertex lies on the left
when crossing $e$ from that face. We denote by $k$ 
this distance and by $m$ the {\it depth} of the other face, say $\mu$  adjacent to $e$ 
in ${\cal M}$. By depth, we mean the number of closing arches (plus possibly the 
number of edges adjacent to the root vertex of ${\cal M}$) separating $\mu$ 
from the external face in the closing process of ${\cal T}$. 
Clearly, a path of length $m$ exists on ${\cal M}$ from the external face to $\mu$.
The fact that $e$ is cut by a leftmost minimal path then imposes 
that $k+1\leq m$, and moreover that $k+1< m$ if the root of ${\cal T}$ lies 
on the side of the white vertex adjacent to $e$ (see Fig.\thtwo). Note finally that, since we 
have chosen an edge with minimal distance, the leftmost minimal path from the external face to $\kappa$ 
cannot cut any edge of ${\cal T}$, thus $k$ is also the depth of $\kappa$ 
(as defined above) in the closing process of ${\cal T}$.
This allows to relate $k$ and $m$ to the charges of the two pieces of ${\cal T}$
on both sides of $e$, namely $k-m=q_\circ$ if the root of ${\cal T}$ 
lies on the side of the black endpoint of $e$ or $m-k=q_\bullet$ if
it lies on the side of the white endpoint of $e$ (see Fig.\thtwo). The above inequalities
on $k$ and $m$ translate into $q_\circ \leq -1$, or equivalently $q_\bullet \geq 2$
(see Fig.\thtwo), which implies that $e$ is non-regular, hence contradicting (S5). 
We deduce that normal edges of ${\cal T}$ are not cut, hence ${\cal T}={\cal T}'$ 
as announced.

To conclude this section, let us note that a direct outcome of the above proof
are the following statements: 

\noindent{\underbar {Proposition (C2):}} 
{\sl Starting from a blossom tree with regular and non regular edges, if we wish
to recover the same tree in a closing-cutting procedure, it is
necessary and sufficient to impose that all the non-regular edges of the tree be marked
as non-sectile. Indeed, all regular edges will remain uncut, even if their
cutting is {\rm a priori} allowed.}
\medskip
\noindent{\underbar {Proposition (C3):}} 
{\sl If the cutting of a map ${\cal M}$ with special edges creates a tree ${\cal T}$ 
in which some of the special edges are regular, then in each cutting of ${\cal M}$
for which some of these edges are no longer marked as special on ${\cal M}$ (i.e.\ become {\rm
a priori} sectile), the latter remain uncut, hence we obtain the {\rm same} tree ${\cal T}$, but
with fewer regular special edges.}

\newsec{Application to bicubic maps with hard particles: proof of main theorem} 

We are now ready to prove our main result \mainth.

\subsec{Admissible maps and equivalence classes of admissible trees}
\fig{An example of two admissible trees (a) and (b) in the same equivalence class, 
i.e.\ whose closing builds the same admissible map (c). Note that the closing of (a) has created
on the map (c) one more NHP edge than that originally on the tree.}{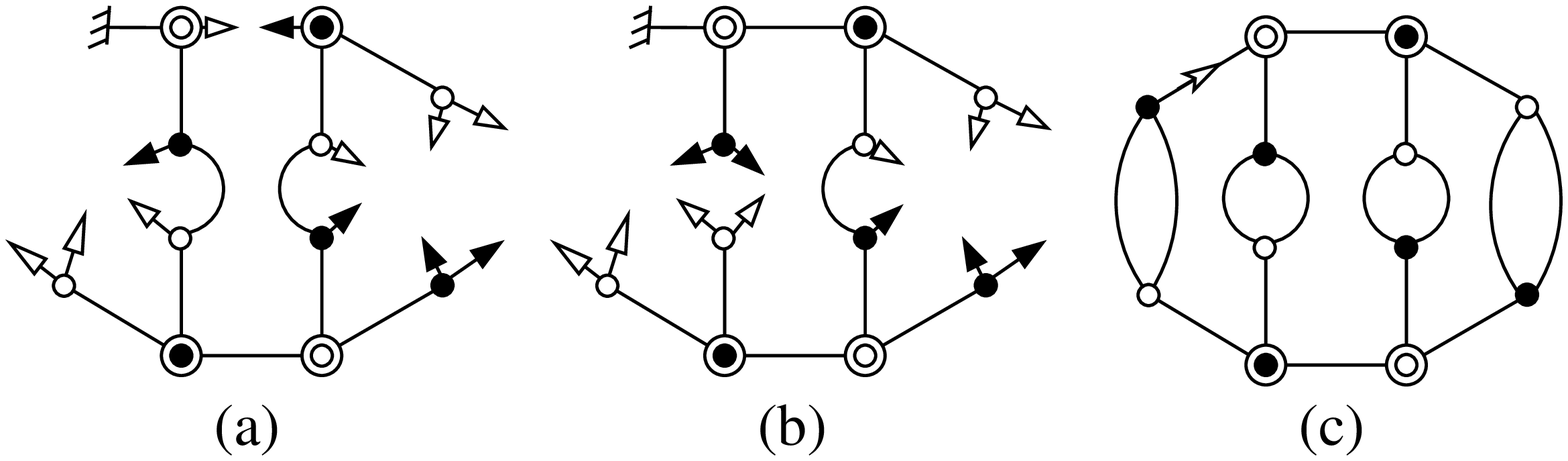}{12.cm}
\figlabel\twoadm
We start with our family of admissible trees ${\cal T}$ as defined by (1)-(6). In particular, 
the regular edges match exactly the HP ones, while the non-regular edges match exactly the 
doubly-occupied ones. We define an {\it admissible map} as any rooted planar bicubic map
${\cal M}$ with particles obtained by the closing of an admissible tree ${\cal T}$. 
A given admissible map may in general be obtained by the closing of several
admissible trees (see Fig.\twoadm). This leads us naturally to define an {\it equivalence relation}
among admissible trees, by ${\cal T} \sim {\cal T}'$ if and only if they lead to the
same map ${\cal M}$. The corresponding equivalence class is denoted by ${\cal C}_{\cal M}$. 
We shall now rely on the propositions (C1)-(C3) above to characterize equivalence classes. 

Starting from an admissible tree ${\cal T}\in {\cal C}_{\cal M}$, we may decide to mark as 
special edges on ${\cal T}$ {\it all its non-regular NHP edges}, keeping all regular HP 
edges as normal edges. The tree ${\cal T}$, together with these special edges, satisfies 
the property (S5). Moreover, the marking of special edges of ${\cal T}$ induces a marking
of special edges on ${\cal M}$ which are necessarily NHP on ${\cal M}$. Note that ${\cal M}$ 
has in general more NHP edges than those of ${\cal T}$, as doubly-occupied edges may be created 
in the closing process (see Fig.\twoadm). Still, according to
proposition (C2), the tree ${\cal T}$ is recovered 
from ${\cal M}$ by the cutting process avoiding the special edges. We immediately deduce 
that any tree in ${\cal C}_{\cal M}$ is obtained from a particular cutting of ${\cal M}$ in which  
a number of its NHP edges have been marked as special. Sweeping all possible
such marking of NHP edges on ${\cal M}$, we get of course in general many more trees than those of 
${\cal C}_{\cal M}$ as some choices of NHP edges on ${\cal M}$ may lead to trees which 
are not admissible. 

\subsec{Cutting of admissible maps} 

Let us now discuss precisely the cutting of an admissible map ${\cal M}$. Denoting 
by ${\rm NHP}({\cal M})$ its set of NHP edges with cardinality
$m\equiv m({\cal M})=|{\rm NHP}({\cal M})|$, we are led to consider the $2^m$ different 
cuttings corresponding to the $2^m$ possible choices of special edges among 
the doubly-occupied ones. Each cutting amounts to the choice of a subset ${\cal S}\subset {\rm NHP}({\cal M})$
with cardinality $p\equiv p({\cal S})=|{\cal S}|$. For any such choice ${\cal S}$, 
the cutting process results in a tree, i.e.\ the marking is acceptable.
Indeed, a non-acceptable choice would require to have a loop made, counterclockwise, of 
either special, hence doubly-occupied edges or white-to-black edges. As black and white
alternate along a loop, we would deduce the presence of a loop made of occupied particles
only. This is clearly impossible as, from the environment of a NHP edge on admissible
trees (see Fig.\envNHP), a NHP edge is always followed or preceded by a HP edge on
the side of its doubly occupied black vertex. Note also that, as it should, the special
edges cannot be adjacent to the root vertex of the map as the latter is unoccupied.
Moreover, the obtained tree ${\cal T}$ satisfies 
the properties (S1)-(S5), from which we immediately deduce the properties (1)-(5) of admissible trees. 
In particular, if an edge is non-regular, it must be special according to (S5), hence it must be 
doubly-occupied. In other words, HP edges on ${\cal T}$ are necessarily regular, which is nothing 
but (5) or equivalently (5'). In general however, the property (6) of admissible trees is 
not guaranteed, i.e.\ the tree may have a number of NHP regular edges. Let us now
discuss the various situations which may occur:
\item{(i)} $m=0$. Then the map satisfies the hard particle constraint. The corresponding
equivalence class is made of the unique tree ${\cal T}$ obtained by cutting ${\cal M}$ with
no special edges. This tree ${\cal T}$ is a good tree. 
\item{(ii)} $m>0$. Then the map does not satisfy the hard particle constraint. The corresponding
equivalence class is made of admissible trees which are all obtained by cuttings of
${\cal M}$ with a subset of its NHP edges marked as special. 
\par 
\noindent In this case (ii), and for a given choice ${\cal S}$ among the $2^m$ possible subsets
of ${\rm NHP}({\cal M})$,  let us further denote by $r\equiv r({\cal T})$ the 
number of NHP regular edges of the resulting tree ${\cal T}$ and by $n\equiv n({\cal T})$ 
its number of NHP non-regular ones. 
We clearly have $n\leq p$ as non-regular edges are necessarily special, and also $p\leq n+r$ 
as all special edges are NHP edges of ${\cal T}$.
Now again two situations may occur:
\item{(ii-a)} $r=0$. Then all NHP edges are non-regular, i.e.\ property (6) is satisfied and
the tree ${\cal T}$ is admissible. In this case, the above inequalities moreover imply that $n=p$,
i.e.\ the NHP edges of ${\cal T}$ are exactly the marked NHP edges of ${\cal M}$.
Note that, from proposition (C1), the marking ${\cal S}$ at hand then corresponds to the {\it only} choice 
(among the $2^m$) that leads to the admissible tree ${\cal T}$.
\item{(ii-b)} $r>0$. In this case, the tree ${\cal T}$ is not admissible. From proposition (C3), we 
also know that the $2^r$ choices corresponding to:
\itemitem{-} marking as special on ${\cal M}$ all the NHP non-regular edges of ${\cal T}$,
\itemitem{-} marking or not the NHP regular edges,
\item{}all lead to the same tree ${\cal T}$. From proposition (C1), these $2^r$ markings are 
the {\it only} choices (among the $2^m$) that lead to ${\cal T}$.
\par
\fig{The $2^m=4$ cuttings of the admissible map ${\cal M}$ in (a) with $m=2$
NHP edges. Three different trees are obtained. A non-admissible tree with $r=1$ 
is obtained for the $2^r=2$ cuttings (b) and (c) corresponding to marking or not the upper
horizontal doubly-occupied edge. The two other possible markings where the lower horizontal
doubly-occupied edge is marked and the upper horizontal doubly-occupied edge is marked
or not lead to two different admissible trees (e) and (f). The equivalent class ${\cal C}_{\cal M}$
is made of these two trees. As one of them has one NHP inner edge and the other two, their 
contributions with a weight $-1$ per NHP inner edge add up to $0$.}{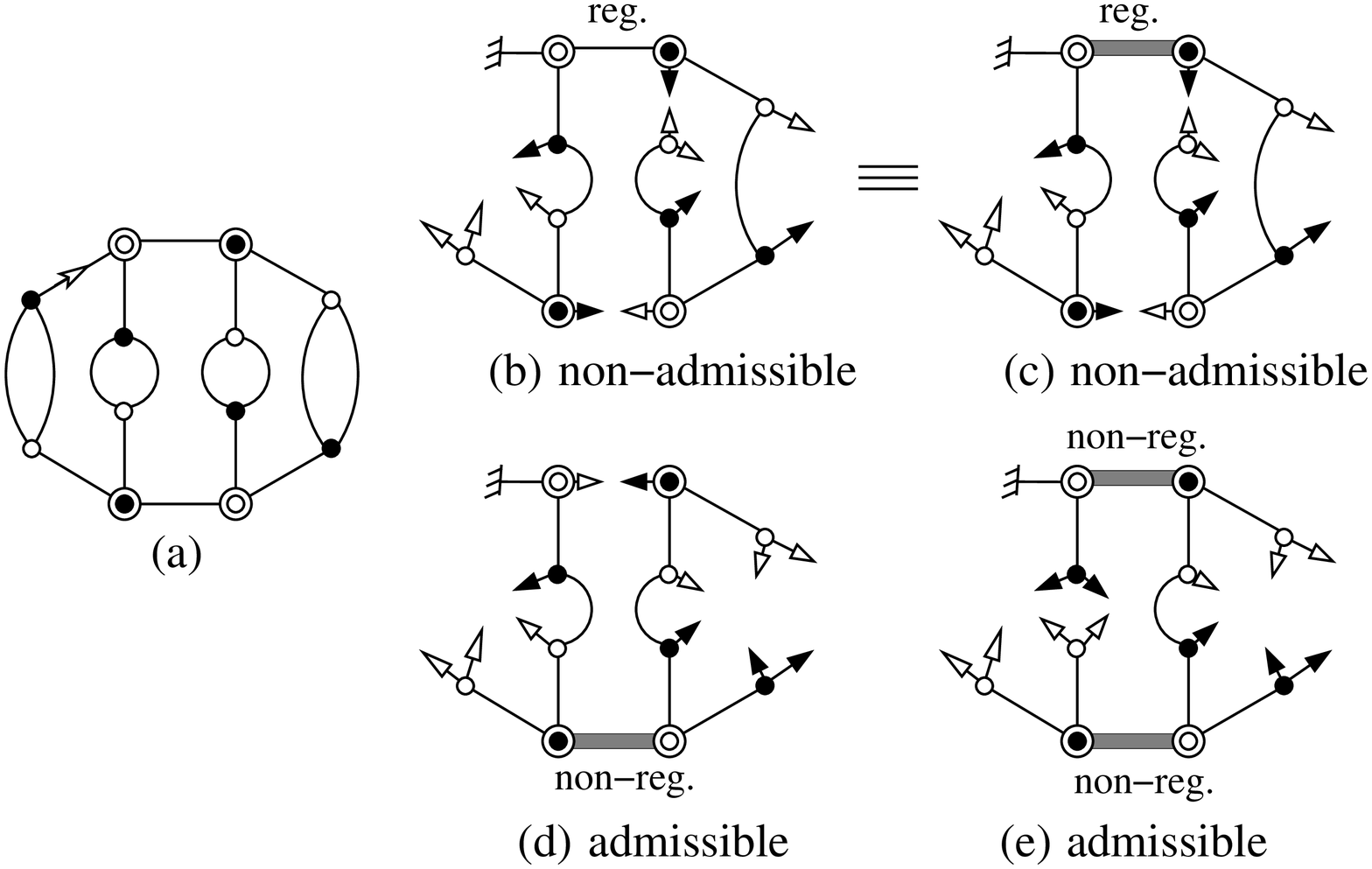}{12.cm}
\figlabel\countex
\fig{The $2^m=4$ cuttings of another admissible map ${\cal M}$ in (a) with $m=2$
NHP edges. Again the class ${\cal C}_{\cal M}$ is made of two trees (b) and (c) whose
contribution with a weight $-1$ per NHP edge add up to $0$. The tree in (d) is identical to
that in (e). It has $r=1$ and is obtained by the $2^r=2$ markings of the original map where
the vertical doubly-occupied edge is selected and the horizontal one is either marked
or not.  As opposed to the case of 
Fig.\countex, the cutting with no marking (b) leads to an admissible tree.}{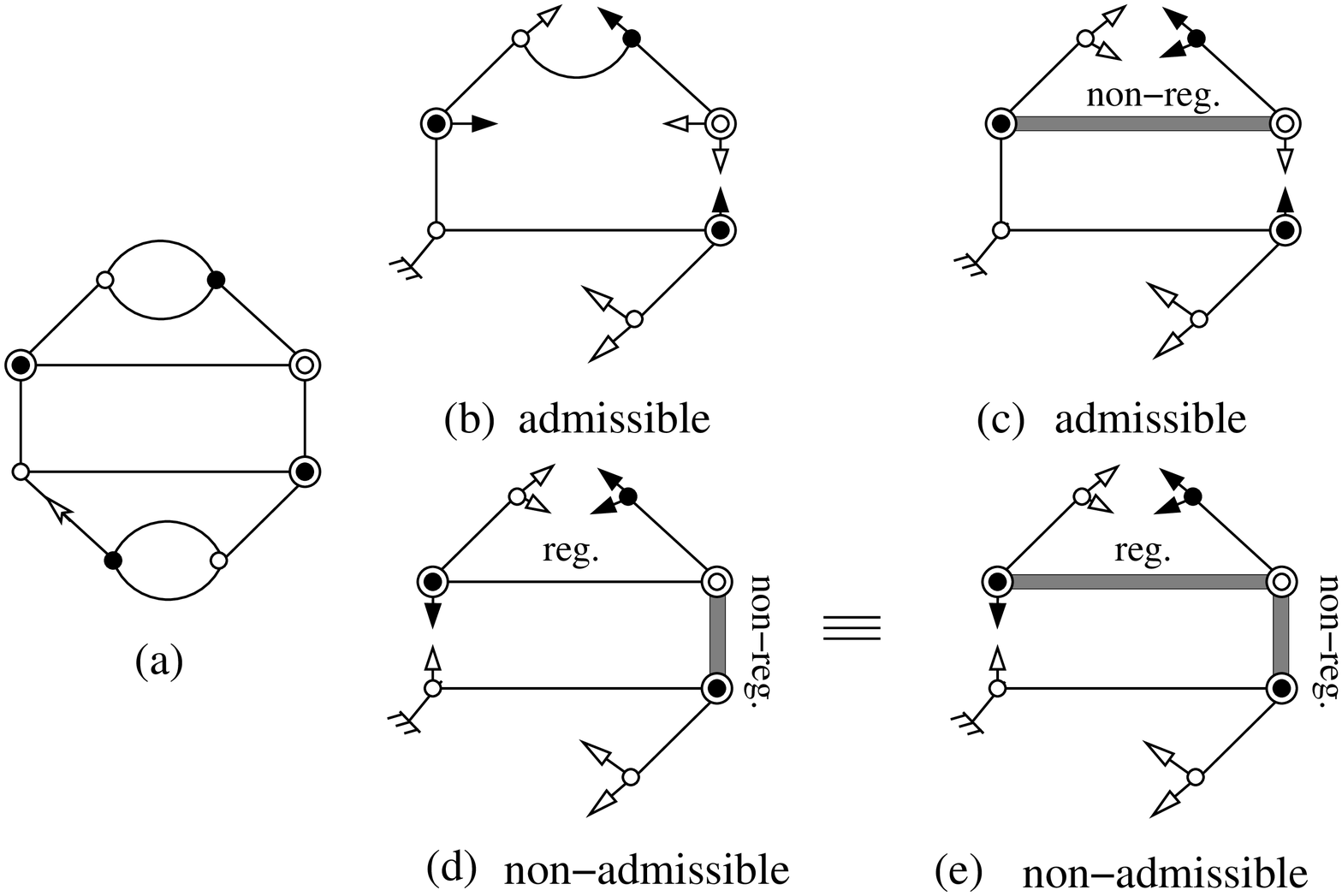}{11.cm}
\figlabel\countextwo
\noindent Fig.\countex\ shows an example of the $2^m=4$ cuttings of an admissible map
with $m=2$ doubly-occupied edges. The resulting trees are two admissible trees, each obtained
for a unique choice of special NHP edges on ${\cal M}$, and a non-admissible
tree with $r=1$ obtained for $2^r=2$ choices of special NHP edges on ${\cal M}$.
The class ${\cal C}_{\cal M}$ is thus made of two elements. Another similar example 
is displayed in Fig.\countextwo.

\subsec{Sum rule and proof of main theorem}

We may now derive the sum rule \mainth\ as follows. We start with
an admissible map ${\cal M}$ with $n_\circ$ white-empty particles, 
$n_\bullet+1$ black-empty particles, $n_\smallcirccirc$ white-occupied particles
and $n_\smallcircbullet$ black-occupied particles. 
We then sum over the $2^{m({\cal M})}$ choices ${\cal S}$ of markings of NHP edges as special edges. 
We assign to each choice a weight $(-1)^{p({\cal S})}$, i.e.\ a weight $-1$ per marking. 
We have clearly:
\eqn\sumone{\eqalign{\sum_{{\cal S}\subset {\rm NHP}({\cal M})}(-1)^{p({\cal S})}& =
\delta_{m({\cal M}),0} \times 1 +(1-\delta_{m({\cal M}),0}) (1+(-1))^{m({\cal M})} \cr
& = \delta_{m({\cal M}),0} \cr}}
as when $m({\cal M})=0$, ${\cal S}={\rm NHP}({\cal M})=\emptyset$ and, when $m({\cal M})>0$, 
each NHP edge may independently be marked or not.
Alternatively we may reorganize the same sum by regrouping all subsets ${\cal S}$
that lead upon cutting to the same tree ${\cal T}$ (admissible or not) and write
\eqn\sumtwo{\eqalign{ \sum_{{\cal S}\subset {\rm NHP}({\cal M})}(-1)^{p({\cal S})}& =
\sum_{\cal T} \sum_{{\cal S}\ {\rm leading} \atop {\rm to\ }{\cal T}}(-1)^{p({\cal S})} 
\cr &= 
\sum_{{\cal T}\ {\rm s.t.} \atop r({\cal T})=0} (-1)^{n({\cal T})} +
\sum_{{\cal T}\ {\rm s.t.} \atop r({\cal T})>0} (-1)^{n({\cal T})} (1+(-1))^{r({\cal T})} \cr
&= \sum_{{\cal T}\in {\cal C}_{\cal M}}  (-1)^{n({\cal T})}
\cr}}
where the sums in the first two lines extend over all trees ${\cal T}$ resulting from all cuttings
of ${\cal M}$. These trees satisfy properties (1)-(5) of admissible
trees, but do not satisfy (6) in general. We have also used the fact that admissible trees 
are those with $r({\cal T})=0$ and the fact that, in this case, $n({\cal T})= p({\cal S})$ 
(cases (i) or (ii)-a above) while trees having $r({\cal T})>0$
are obtained in $2^{r({\cal T})}$ ways corresponding to marking all their non-regular NHP
edges and marking or not independently each regular NHP edge. Comparing Eqs. \sumone\ and
\sumtwo\ leads directly to the sum rule 
\eqn\refined{\sum_{{\cal T}\in {\cal C}_{\cal M}}(-1)^{n({\cal T})}= \delta_{m({\cal M}),0}}
As $m({\cal M})=0$ precisely means that the map satisfies the hard particle constraint, 
the r.h.s. above is nothing but the characteristic function of rooted bicubic maps with hard particles
within the larger set of admissible maps. Summing over all admissible maps ${\cal M}$, 
we immediately deduce the main theorem \mainth.

\newsec{Discussion}

In this paper, we have presented a direct combinatorial derivation of the generating
function for hard particles on planar bicubic maps by use of admissible trees. 
To conclude, we shall first briefly present yet another application of our construction 
to the Ising model on tetravalent maps, described by a similar chain-interacting four-matrix integral. 
More generally, our construction may be extended to models of maps with
particles subject to exclusion rules, all described, as mentioned in the introduction,
by chain-interacting multi-matrix integrals. We shall sketch below the general 
structure of admissible trees for the general multi-matrix case. We will finally comment
on the predicted existence of an Ising-like crystallization transition for a particular 
model of trees with particles, with possible re-interpretation as a branching process.

\subsec{The Ising model revisited}

The Ising model on tetravalent maps is well known to be described by a two-matrix integral
with quartic potential. For the present purpose however, it is more instructive to
view it as a hard particle model on bipartite maps, hence described by a four-matrix
integral. More precisely, let us consider hard particles on bipartite planar maps
with bi- and tetravalent vertices only and with the constraint that all tetravalent
vertices are unoccupied, while all bivalent ones are occupied. In the matrix language,
the corresponding action reads
\eqn\fourmatising{U(M_1,\ldots,M_4)=M_1 M_2-M_2 M_3+M_3 M_4 -g \left({M_1^4\over 4}+{M_4^4\over 4}\right)\
-g\,z \left( {M_2^2\over 2}+{M_3^2\over 2}\right)}
with four matrices describing black-empty tetravalent ($M_1$), white-occupied bivalent ($M_2$), 
black-occupied bivalent ($M_3$) and white-empty tetravalent ($M_4$) vertices and with
a weight $g$ per vertex and $z$ per particle. The Ising model is recovered by erasing
all bivalent vertices. This results in tetravalent maps with black and white vertices,
all with weight $g$, linked irrespectively of their color, but with an effective edge coupling $1$ 
for edges adjacent to vertices of opposite color, and $g z$ for edges adjacent to
vertices of the same color. This is nothing but the Ising model upon interpreting
colors as spins.

Returning to our bipartite hard-particle formulation, it proves convenient to consider maps 
with a single additional bivalent black-empty root vertex. This corresponds to what
was called ``quasi-tetravalent" maps in Ref.\BMS. The corresponding maps may be
enumerated along the same lines as above by counting properly signed admissible trees.
An admissible tree is now defined as follows:
\item{(1)} Its inner vertices are of four types: black or white vertices, each of which being
either empty or occupied by a particle. All empty vertices have valence four and all
occupied ones have valence two.
\item{(2)} It is bicolored, i.e.\ all its inner edges connect only black (empty or occupied)
inner vertices to white (empty or occupied) ones.
\item{(3)} It has two kinds of leaves: the so-called buds, connected only to black (empty or occupied)
inner vertices and the so-called leaves, connected only to white (empty or occupied) inner vertices.
We attach a {\it charge} $-1$ to each bud and a charge $+1$ to each leaf.
\item{(4)} Its total charge ($\#$leaves $-$ $\#$buds) is equal to $2$.
\item{(5)} All HP edges of the tree are regular.
\item{(6)} All NHP edges of the tree are non-regular.
\par 
\noindent As a consequence of (1)-(4), any inner edge separates the tree into a piece of total charge 
$q_\circ\equiv 1{\rm \ mod\ }2$ starting with a white vertex, and a piece of total charge 
$q_\bullet=2-q_\circ\equiv 1{\rm \ mod\ }2$ starting with a black vertex. Such an edge will be again called 
of type $(q_\circ:q_\bullet)$. As before, the last two characterizations (5) and (6) may be rephrased into:
\item{(5')} Every HP edge is of type $(q_\circ:q_\bullet)$ with $q_\circ\geq 1$,
hence $q_\bullet \leq 1$.
\item{(6')} Every NHP edge is of type
$(q_\circ:q_\bullet)$ with $q_\circ\leq -1$, hence $q_\bullet \geq 3$.
\par
\noindent A detailed survey of possible local environments shows that admissible trees have only
NHP edges of type $(q_\circ=-1:q_\bullet=3)$, while HP edges are of type
$(q_\circ=3:q_\bullet=-1)$ or $(q_\circ=1:q_\bullet=1)$. As before, the generating function
$G_{\rm QTIsing}$ of quasi-tetravalent maps with Ising spins
is expressed as that of admissible trees with a weight $-1$ per NHP edge. The latter
is easily obtained from the generating function of planted trees, expressed in terms 
of generating functions for half trees. The half trees to be considered are 
\item{(i)} half-trees of charge $-1$ attached to a white-occupied vertex, generated by $A$
\item{(ii)} half-trees of charge $1$ attached to a white-occupied vertex, generated by $B$
\item{(iii)} half-trees of charge $3$ attached to a white-occupied vertex, generated by $C$
\item{(iv)}  half-trees of charge $1$ attached to a black-empty vertex, generated by $R$
\item{(v)}  half-trees of charge $3$ attached to a black-empty vertex, generated by $S$
\item{(vi)} half-trees of charge $-1$ attached to a white-empty vertex, generated by $U$
\item{(vii)} half-trees of charge $1$ attached to a white-empty vertex, generated by $V$
\par
\fig{A pictorial representation of the recursive relations (6.2) for half-trees, 
respectively attached to a white-occupied vertex, with charge $-1$ ($A$), $1$ ($B$) or
$3$ ($C$), to a black-empty vertex, with charge $1$ ($R$) or $3$ ($S$) and to a 
white-empty vertex, with charge $-1$ ($U$) or $1$ ($V$). The various terms collect 
all possible environments compatible with the charge characterization of admissible trees. 
We also indicate the various weights $g$, $z$ and minus signs.}{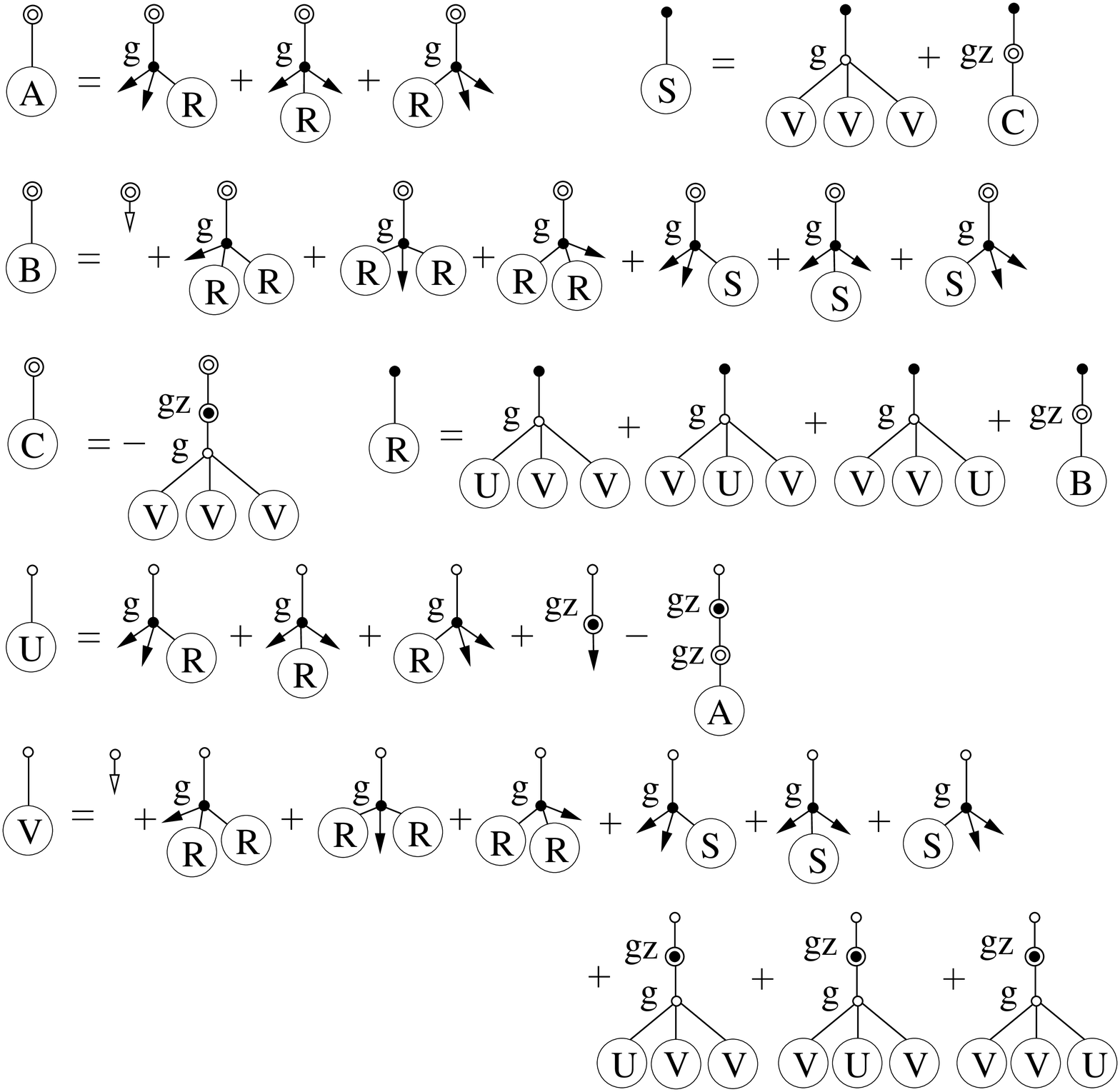}{13.cm}
\figlabel\recurising
\noindent These functions satisfy recursion relations, obtained by inspection of the environment of
the root, as displayed in Fig.\recurising:
\eqn\isingrecur{\eqalign{
A&= 3 g\, R \cr B&= 1 + 3 g\, R^2 + 3 g\, S \cr C &= -g^2z\, V^3 \cr R &= 3 g\, V^2 U+gz\, B\cr
S&= g\,V^3+gz\, C\cr U&= 3 g\, R+gz\,-g^2z^2\, A\cr V&= 1+3 g\, R^2 + 3 g\, S + 3 g^2z\, V^2 U\cr}}
The generating function $G_{\rm QTIsing}$ is easily related to the above, now
via 
\eqn\GQTI{G_{\rm QTIsing}(g,z)=g R- g^2 \left(R^3+6 R S + g z\, V^3\right)}
obtained by attaching admissible trees to either a leaf (first contribution) or a bud
(second contribution with minus sign).
To make the connection with previous results on the Ising model, we note that, by
elimination, the system \isingrecur\ reduces to a single algebraic equation for $V$, which,
upon the substitutions $P=V/(1-g^2 z^2)$, $x=g(1-g^2 z^2)^2$ and $v=g z$, reads
\eqn\eqforP{P=1+3 x^2 P^3 +v^2 {P\over (1-3 x P)^2}}
This equation is nothing but that of Refs.\ISING\ and \BMS. We also easily check that
the formula \GQTI\ matches that of Ref.\BMS.

\subsec{Generalizations}

The construction of this paper may be extended so as to describe hard
particles on bipartite maps with vertices of arbitrary valences (arbitrary polynomial 
potential for each matrix in the action). In this case, the admissible trees 
themselves have arbitrary valences. Again their fundamental characterization is that
HP edges have to be regular while NHP ones have to be non-regular. A weight $-1$
per NHP edge is still required to perform the right counting.

More interestingly, the method can easily be extended to describe particles with
exclusion rules on bipartite maps. As discussed in the introduction, 
we may impose that each edge have a total number of adjacent particles of
at most $k$, corresponding to a multi-matrix integral over $2k+2$ matrices.
Very generally, the tree formulation will use admissible trees with vertices occupied
by at most $k$ particles, and whose edges either satisfy the exclusion rule (ER edges) or
do not (NER edges). The admissible trees are now required to have all their ER edges
regular and all their NER ones non-regular. For maps with bounded valences, the number 
of allowed $(q_\circ:q_\bullet)$ pairs will remain finite, hence the number of
half-trees to be dealt with is finite as well. To recover the correct generating
function for the maps, we simply have to include again a weight $-1$ per NER edge.
In practice, the enumeration might prove cumbersome, although straightforward.
An example with $k=2$ of a six-matrix model which displays the tricritical Ising 
transition was solved in Ref.\CRI\ by matrix techniques. The corresponding 
algebraic equations may be easily re-interpreted in terms of generating functions for
(signed) admissible half-trees.
The case on non-bipartite maps is more subtle and so far, we only know how to treat
the case of odd $k=2p+1$ upon transforming its configurations into those 
of bipartite maps with edges adjacent to at most $p$ particles.

\subsec{Trees with Ising transition}

The generating function $G_{\rm BMHP}$ of planar bicubic maps with hard particles, rooted
at a black-empty vertex is also that of good trees (up to a multiplicative factor $g$ for
the removed root vertex). Recall that we may view the good trees 
as those obeying: 
\item{(G1)} the hard particle constraint is satisfied on the tree,
\item{(G2)} the properties (1) to (5) of Sects.3.1 and 3.2 hold,
\item{(G3)} the closing of the tree produces a map satisfying the hard particle constraint.
\par 
\noindent These trees are moreover rooted
at one of their three unmatched leaves. We may as well consider the same trees, but now planted
at {\it any} of their leaves, and generated by $G_{\rm leaf}$. From the
relation $G_{\rm leaf}=R$, apparent from Figs.\figrecur\ and \planted, we can interpret $R$ 
as generating {\it planted good trees} satisfying properties (G1)-(G3), with weights
$g$ per vertex and $z$ per particle. In particular, as we restrict ourselves to good trees 
{\it ab initio}, no sign factors are required. We have the relation
\eqn\conjug{G_{\rm BMHP}\vert_{g^{2n}}={3\over n+2} R\vert_{g^{2n-1}}}
between the coefficients in $g$ of the two series, expressing that only $3$ of the $n+2$ 
available leaves are unmatched in good trees with $2n-1$ vertices. Note also that the
property (5) may be rephrased into: all descending subtrees of a planted good tree have
charge $q_\bullet=1$ if they start by a black vertex, and charge
$q_\circ=2$ if they start by a white vertex.

Alternatively, we may also interpret $S$ as counting a slightly different version
of planted good trees obtained as follows. Starting from a planar bicubic
map with hard particles rooted at a black-empty vertex, and using the cutting 
process of Fig.\decoupagetwo\ based on the geodesic tree, we may stop the process
at step (b), i.e. keep the root vertex and its two attached buds and plant the
tree at one of these buds.
The resulting tree is a new type of planted good tree satisfying the following properties
(G1')-(G3'). The property (G1') is the same as (G1). The property (G2') is similar to
(G2) but with (4) replaced by: its total charge is $0$ and with (5) replaced by: 
all descending subtrees have charge $q_\bullet=1$ if they start by a 
black vertex, and charge $q_\circ=2$ if they start by a white vertex. Finally (G3') states
that closing the tree into a map by simply connecting iteratively all successive bud-leaf pairs
produces a map with no NHP edge. It is then easy to see that these planted trees, with weights
$g$ per vertex and $z$ per particle are generated by $S$, as apparent from Fig.\figrecur\ by
replacing the white occupied root vertex of $S$ by a bud (for convenience, we also
include an empty tree with weight $1$). This now results in the relation
\eqn\conjug{G_{\rm BMHP}\vert_{g^{2n}}={3\over 2(n+2)} S\vert_{g^{2n}}}
for $n>0$, with an extra factor of $1/2$ to account for the two choices of bud where to plant
the tree.

For $g$ and $z$ non negative, the good trees generated by $S$ form a family of planted trees 
with {\it local non negative weights} only, characterized by a number of local constraints, 
and subject to a non-trivial {\it global} ``goodness" constraint (G3') that, upon closing, 
no NHP edge is created.
The singularity structure of $S$ is inherited from that of $G_{\rm BMHP}$, hence the statistics
of these good trees undergoes at $z=z_{+}$ the Ising-like crystallization transition
of Sect.2.3. 

As the weights are non-negative, we may rephrase them in terms of probabilities by assigning
to each such planted good tree ${\cal T}$ with $2n({\cal T})$ vertices and $h({\cal T})$ particles
a probability 
\eqn\probaT{p({\cal T})={1\over S} g^{2n({\cal T})} \, z^{h({\cal T})}}
By construction, these add up to $1$ on the set of planted good trees. 
\fig{Local evolution rules for the branching process associated to good trees 
without particles (see text).}{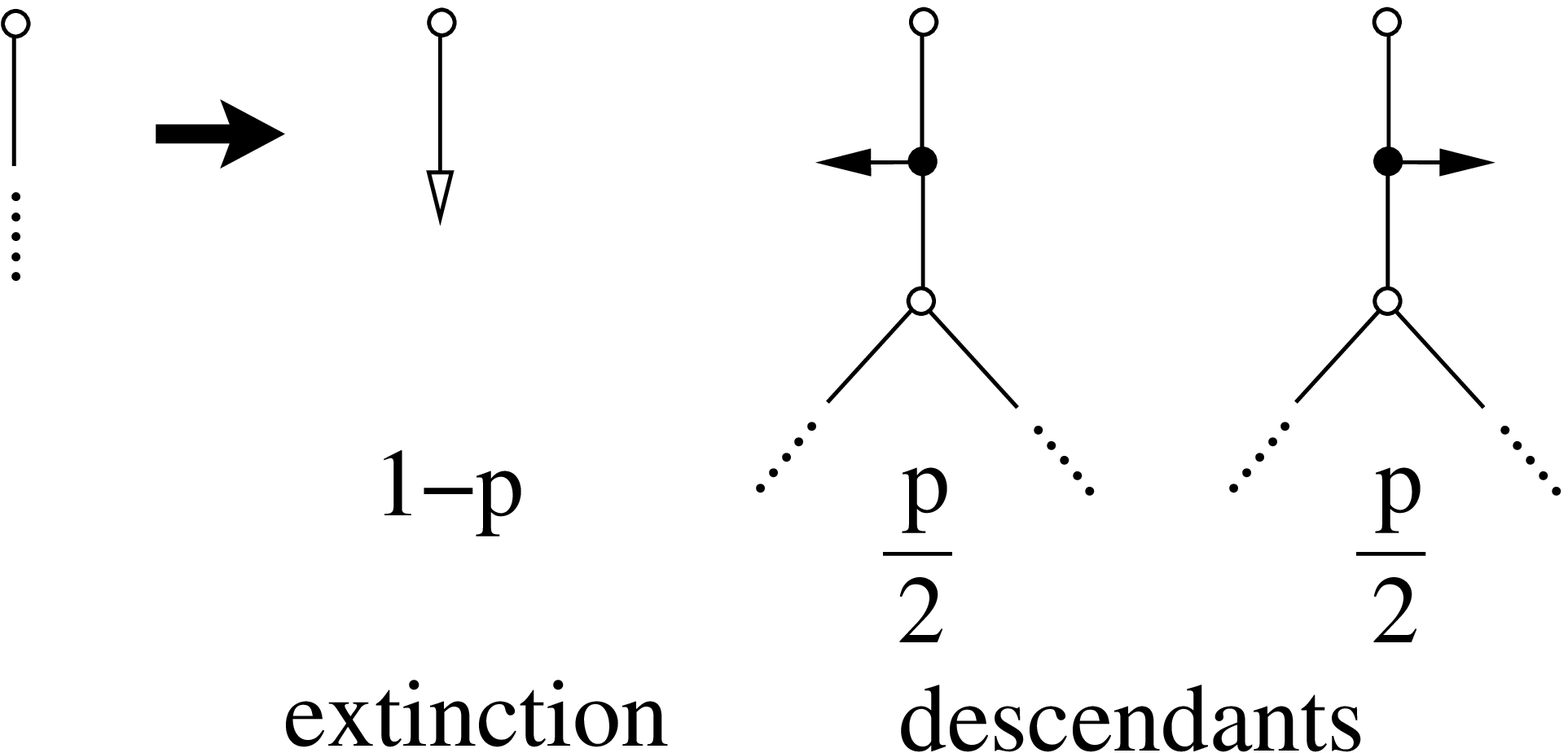}{7.cm}
\figlabel\evol
In the case $z=0$ of maps without particles, the planted trees are generated by $S_0\equiv S(z=0)$,
solution of $S_0=1+2 g^2 S_0^2$, namely $S_0=(1-\sqrt{1-8 g^2})/(4 g^2)$, with a singularity
at $g^2=1/8$. Their probabilities may be associated to those of a branching process as follows.
Each tree is interpreted as the genealogical tree of a species whose individuals
live at the level of the edges descending from white vertices and are subject to the
evolution rules sketched in Fig.\evol. The individuals either die without
descendants with probability $(1-p)$ (if the edge at hand is a leaf) or, with probability
$p$, first make a two-state choice with probability $1/2$ (black vertex 
with the bud pointing to the left or to the right) and then give rise to two children (the
two descending edges of the white descendant vertex of the black one).
Comparing with the recursive rules for good trees without particles, this amounts to 
setting $g^2=p(1-p)/2$. This leads to $S_0=1/(1-p)$ if we assume $0<p<1/2$, and the 
probability \probaT\ of a tree ${\cal T}$ reads:
\eqn\probaB{p({\cal T})=\left( {1\over 2}\right)^{n_\bullet({\cal T})} p^{n_\circ({\cal T})}
(1-p)^{n_{\rm leaf}({\cal T})}}
if the tree has $n_{\rm leaf}$ leaves, $n_\bullet$ black vertices and $n_\circ$ white vertices.
This probability matches clearly that of the branching process. 

Note that for $p>1/2$, the formula \probaB\ still corresponds to the
probability of the branching process, but it does not match the probability \probaT\
for the tree as $S_0$ is now given by $1/p$. The relative factor $(1-p)/p$ is nothing but
the probability that the branching process remain finite, i.e. that the species become extinct.
It is known that this probability is less than $1$ for $p>1/2$. This is a general feature
of so-called Galton Watson branching processes, with local evolution rules \GALWA. These processes
undergo a transition at the value of the parameters (here $p$) at which the average number
of children (here equal to $2p$) attains $1$. Note that the transition takes place precisely
at the singularity of $S_0$ as $g^2=1/8$ when $p=1/2$. This should not come as a surprise as,
when $g^2\to (1/8)^{-}$, the average size of the trees diverges, which is the signal that
infinite branching processes are about to appear. Note that the probabilities \probaT, 
which add up to $1$ by construction, may still be interpreted as those of the branching
process provided the latter is {\it conditioned} to remain finite (i.e. get extinct).

With only probabilistic weights and local evolution rules, branching processes are trapped
in the ``universality class" of the Galton Watson transition. In the tree language, this 
amounts to always having the same (square-root) singularity structure for generating
functions. In the presence of particles, good trees however seem to be a way to escape from this 
Galton Watson framework as we know from the above analysis that they undergo an Ising-like 
crystallization transition at which the singularity has a different nature, namely that of 
a cubic root responsible for the exponent $7/3$ in Eq.\Ggrowthmc.
The new ingredient that gives rise to this new universality class is the non-local 
goodness constraint on the trees. 
We may again view the good trees and their probabilities \probaT\ as describing 
branching processes with local evolution rules, again conditioned to remain finite,
and further conditioned by the 
(somewhat non-natural) global goodness constraint. A good understanding of
this branching process is still lacking. As before, approaching the critical line
$g\to g_c(z)$ must correspond to attaining an average number of children equal to $1$. It
would be interesting to understand:
\item{(i)} whether or not the multi-critical point at $z=z_+$ may be attained in the
physically acceptable range of parameters,
\item{(ii)} if so, what general property of the branching process
governs the approach to the multi-critical point.
\par
\noindent To this end, it seems likely that a formulation of maps as labelled trees like that
described in Refs.[\xref\CS,\xref\MOB] rather than blossom trees may be more appropriate.
Also, a generalization of the goodness condition to infinite trees would be desirable.

On the other hand, the use of the larger class of admissible trees allows to get rid
of the global goodness constraint, but at the expense of having negative weights.
Still we might hope to be able to re-interpret these signs as subtractions of 
overcounted branching processes. In the present blossom tree language however,
the trees in the same equivalence classes of a non-admissible map (which add up to zero),
have very different shapes and we have not been able to rearrange them to give
such a re-interpretation.

\newsec{Conclusion}
In this paper, we have shown how to extend the bijective methods of Refs.[\xref\CORV-\xref\CHP]
to the case of planar bicubic maps with hard particles. In the matrix language, this extends
the bijective combinatorial treatment of the one- and two-matrix models to the case of 
a four-matrix integral with action \fourmat\ displaying a chain-like interaction. 
The first lesson to be drawn from this study is that the bijective approach still
applies in this case although the matrix model solution clearly indicated the presence of
signed weights in the auxiliary counting functions (cf. Eq.\matrecur). We have indeed shown 
that these signs
could be interpreted in terms of a kind of inclusion-exclusion principle, here projecting
the generating function of admissible trees onto that of good trees, hence reproducing
$G_{\rm BMHP}(g,z)$. A second lesson is that the cutting procedure of Ref.\SCHth, based on
paths which are minimal with respect to the canonical geodesic distance on the map, may
be extended for each given map to other (less canonical) definitions of the geodesic
distance. In particular, the introduction of blocked edges was instrumental in our
proof. Physically, this may be understood as a form of interaction between the space metric 
(distances on the map) and matter (presence of particles correlated to the blocked edges).
This freedom in the choice of the definition of distances may prove useful to tackle
other problems, for instance extended hard objects such as hard dimers (occupied edges), 
animals (occupied clusters), etc...

An interesting property of this study is the robustness of the method, which may
be adapted to the case of arbitrary potentials as well as higher order exclusion rules
(higher number of matrices). This extends {\it de facto} the class of problems amenable
to a bijective treatment to a much larger family of models. The latter is expected
to include all possible (multi-)critical points of two-dimensional statistical models
on random surfaces whose universality classes are described by minimal conformal field 
theories coupled to two-dimensional quantum gravity. All these models may therefore
be understood in terms of trees with proper signed weights.
In particular, this combinatorial framework provides a unified formulation of
both the Ising model and the model of hard particles on bipartite graphs, which 
eventually explains their common universality class.

Finally, the existence of a tree formulation implies a connection between these
map enumeration problems and spatial branching processes. We may hope that these
two fields of research may benefit from this connection. In particular, the critical
points of maps with generalized particle exclusion rules lying in unitary conformal
classes should translate into well defined branching processes (with positive
probability weights). By a global conditioning, we might even be able to escape from 
the restricted class of Galton-Watson processes and reach in the continuum limit
new classes of continuous random trees, generalizing that of Ref.\Aldous.

\bigskip
\noindent{\bf Acknowledgments:} 
All the authors acknowledge the support of the 
Geocomp project (ACI Masse de donn\'ees). J.B. acknowledges
financial support from the Dutch Foundation for Fundamental Research on Matter
(FOM). P.D.F. acknowledges support from the European network ``Enigma", grant
MRTN-CT-2004-5652.

\listrefs
\end